\documentclass[12pt]{imsart}
\usepackage{amsthm,amsmath}
\usepackage[colorlinks,citecolor=blue,urlcolor=blue]{hyperref}

\usepackage{amsmath,amsfonts,amssymb,amsthm}

\usepackage{enumitem}

\RequirePackage{amsthm,amsmath,amsfonts,amssymb}
\RequirePackage[numbers]{natbib}
\RequirePackage[colorlinks,citecolor=blue,urlcolor=blue]{hyperref}
\RequirePackage{graphicx}
\usepackage{lineno,hyperref}

\startlocaldefs
\modulolinenumbers[5]

\usepackage[utf8]{inputenc}
\usepackage[english]{babel}

\usepackage[margin=1in]{geometry}
\usepackage{bm}
\usepackage{graphicx}
\usepackage{amssymb,amsmath,amsthm,mathrsfs}
\usepackage{epstopdf}
\usepackage{algorithm}
\usepackage{amsfonts,delarray}
\usepackage{algorithm,algcompatible,amsmath}
\usepackage{rotating,color}
\usepackage{subfigure}
\usepackage{multirow}
\usepackage{titlesec,textcase}
\usepackage{longtable}
\usepackage{bbm}
\usepackage{rotating}

\newtheorem{Theorem}{Theorem}[section]
\newtheorem{Lemma}[Theorem]{Lemma}
\newtheorem{Corollary}[Theorem]{Corollary}

\newtheorem{Definition}[Theorem]{Definition}

\theoremstyle{definition}

\RequirePackage[normalem]{ulem} 

\definecolor{rp}{RGB}{83,54,106}

\def\boxit#1{\vbox{\hrule\hbox{\vrule\kern6pt\vbox{\kern6pt#1\kern6pt}\kern6pt\vrule}\hrule}}

\endlocaldefs

\allowdisplaybreaks

\begin{document}
\begin{frontmatter}
\title{Central limit theorem for the average closure coefficient }

\runtitle{ distribution of the local closure coefficient}
\runauthor{ }
\begin{aug}

\author[A]{\fnms{Mingao} \snm{Yuan}\ead[label=e1]{mingao.yuan@ndsu.edu}}



\address[A]{Department of Statistics,
North Dakota State University,
\printead{e1}}
\end{aug}

\begin{abstract}
Many real-world networks exhibit the phenomenon of edge clustering, which is typically measured by the average clustering coefficient. Recently, an alternative measure, the average closure coefficient, is proposed to quantify local clustering. It is shown that the average closure coefficient possesses a number of useful properties and can capture complementary information 
missed by the classical average clustering coefficient.
In this paper, we study the asymptotic distribution of the average closure coefficient of a heterogeneous  Erd\"{o}s-R\'{e}nyi random graph. We prove that the standardized average closure coefficient converges in distribution to the standard normal distribution. In the Erd\"{o}s-R\'{e}nyi random graph, the variance of the average closure coefficient exhibits the same 
phase transition phenomenon as the average clustering coefficient.
\end{abstract}

\begin{keyword}[class=MSC2020]
\kwd[]{60K35}
\kwd[;  ]{05C80}
\end{keyword}

\begin{keyword}
\kwd{average closure coefficient}
\kwd{random graph}
\kwd{asymptotic distribution}
\end{keyword}

\end{frontmatter}

\section{Introduction}
\label{S:1}

A network or graph $G=(V,E)$ is a pair of node set $V$ and edge set $E$. The edges in $E$ represent the interactions between nodes.
Networks are widely used to understand and model many complex systems \cite{CGL16,N03}.
In sociology, relationships among social actors can be depicted by networks and network analsyis is used to study structures of interdependencies among social units \cite{BM08}. In biology, network is an important tool for understanding how the interactions between molecules determine the function of cells \cite{BM03}. In psychology, network analysis is applied to identify and analyse patterns of pairwise conditional
dependencies in multivariate psychological data \cite{BM01}.

Many real-world networks exhibit a trait that the edges tend to cluster. For instance, in a social network, the friends of
a friend are more likely to be friends \cite{BW00}.  In co-authorship
networks, the collaborators of an author tend to be co-authors \cite{JGN01}.  The average clustering coefficient is the standard metric to measure the extent of clustering \cite{WS98}. The local clustering coefficient of a node is defined as the fraction of pairs of
its neighbors that are connected by an edge. The average clustering coefficient is the average of the local clustering coefficients of all nodes. The average clustering coefficient has wide applications 
in network data analysis \cite{BBA01,CLX01,TPL01}.

Recently, \cite{YBL19} introduces an alternative metric, the average closure coefficient, to measure the extent of clustering of a network. The local closure coefficient of a node is the fraction of closed wedges that emanate from the node. The average of the local closure coefficients of all nodes is called the average closure coefficient. It is shown that the local (or average) closure coefficient has remarkably different properties from the local (or average) clustering coefficient \cite{YBL19}. For example, as the degree of a node increases, the local closure coefficient tends
to increase, but the average clustering coefficient tends
to decrease \cite{YBL19}. Therefore, the closure coefficients can capture complementary information missed by the clustering coefficients and are
useful in link prediction,
role discovery, outlier detection, etc.

Understanding the asymptotic properties of network statistics is a fundamental research topic in network analysis \cite{BM06,CHHS20,Y23,Y23b,Y23c, ZY23,YZ23}. In this paper, we are interested in the limiting distribution of the average closure coefficient in a heterogeneous  Erd\"{o}s-R\'{e}nyi random graph. The average closure coefficient is a sum of dependent terms. The classic central limit theorem can not be directly applied to obtain its asymptotic distribution. We prove that the standardized average closure coefficient converges in distribution to the standard normal distribution. We find that the variance of the average closure coefficient  exhibits the same phase transition phenomenon as the  average clustering coefficient.

The rest of the paper is organized as follows. In section 2, we introduce the definition of the average closure coefficient and the heterogeneous Erd\"{o}s-R\'{e}nyi random graph and present the main result. The proof is deferred to section 3.

\medskip

\noindent
{\bf Notation:} We adopt the  Bachmann–Landau notation throughout this paper. Let $a_n$  and $b_n$ be two positive sequences. Denote $a_n=\Theta(b_n)$ if $c_1b_n\leq a_n\leq c_2 b_n$ for some positive constants $c_1,c_2$. Denote  $a_n=\omega(b_n)$ if $\lim_{n\rightarrow\infty}\frac{a_n}{b_n}=\infty$. Denote $a_n=O(b_n)$ if $a_n\leq cb_n$ for some positive constants $c$. Denote $a_n=o(b_n)$ if $\lim_{n\rightarrow\infty}\frac{a_n}{b_n}=0$. Let $\mathcal{N}(0,1)$ be the standard normal distribution and $X_n$ be a sequence of random variables. Then $X_n\Rightarrow\mathcal{N}(0,1)$ means $X_n$ converges in distribution to the standard normal distribution as $n$ goes to infinity. Denote $X_n=O_P(a_n)$ if $\frac{X_n}{a_n}$ is bounded in probability. Denote $X_n=o_P(a_n)$ if $\frac{X_n}{a_n}$ converges to zero in probability as $n$ goes to infinity. Let $\mathbb{E}[X_n]$ and $Var(X_n)$ denote the expectation and variance of a random variable $X_n$ respectively. $\mathbb{P}[E]$ denote the probability of an event $E$. Let $f=f(x)$ be a function. Denote $f^{(k)}(x)=\frac{d^kf}{dx^k}(x)$ for any positive integer $k$. $\exp(x)$ denote the exponential function $e^x$.
For positive integer $n$, denote $[n]=\{1,2,\dots,n\}$. Given a finite set $E$, $|E|$ represents the number of elements in $E$. For positive integers $i,j,k$, $i\neq j\neq k$ means $i\neq j, j\neq k, k\neq i$. Given positive integer $t$, $\sum_{i_1\neq i_2\neq\dots\neq i_t}$ means summation over all integers $i_1,i_2,\dots,i_t$ in $[n]$ such that $|\{i_1,i_2,\dots,i_t\}|=t$. $\sum_{i_1< i_2<\dots< i_t}$ means summation over all integers $i_1,i_2,\dots,i_t$ in $[n]$ such that $i_1<i_2<\dots<i_t$.  For two sets $U,V$, $U-V$ represents the set of elements in $U$ but not in $V$. $U+V$ means the union of $U$ and $V$.

\section{Main results}

An \textit{undirected} graph on $\mathcal{V}=[n]$ is the pair $\mathcal{G}=(\mathcal{V},\mathcal{E})$, where $\mathcal{E}$ is a set of subsets with cardinality 2 of $\mathcal{V}$. The elements in $\mathcal{V}$ are called nodes or vertexes and elements in $\mathcal{E}$ are called edges. A symmetric adjacency matrix $A$ is usually used to represent a graph. In $A$, $A_{ij}=1$ indicates $\{i,j\}$ is an edge. Otherwise,  $A_{ij}=0$. The degree $d_i$ of node $i$ is the number of edges connecting it, that is, $d_i=\sum_{j=1}^nA_{ij}$.  
A wedge or 2-path in a graph is two edges that share exactly
one common node. The common node is called
the center of the wedge. A wedge is closed if nodes at both ends are connected by an edge. For example, edge $\{1,2\}$ and edge $\{2,3\}$ form a wedge. Node 2 is the center of this wedge.  The head of this wedge is node 1 or 3. If  nodes 1 and 3 are connected by an edge, that is, $\{1,3\}$ is an edge, then the wedge is closed.
 A graph is said to be random if $A_{ij} (1\leq i<j\leq n)$ are random.

\begin{Definition}\label{def1}
Let $\alpha$ and $\beta$ be constants between zero and one, that is, $\alpha,\beta\in(0,1)$, and 
\[W=\{w_{ij}\in[\beta,1]| 1\leq i,j\leq n, w_{ji}=w_{ij} ,w_{ii}=0\}.\]
Define a heterogeneous random graph  $\mathcal{G}_n(\alpha,\beta, W)$ as  
\[\mathbb{P}(A_{ij}=1)=p_nw_{ij},\]
where $A_{ij}$ $(1\leq i<j\leq n)$ are independent, $A_{ij}=A_{ji}$ and $p_n=n^{-\alpha}$.
\end{Definition}

In $\mathcal{G}_n(\alpha,\beta, W)$, $\mathbb{E}[d_i]=\sum_{k=1}^np_nw_{ik}$. Generally speaking, $\mathbb{E}[d_i]\neq \mathbb{E}[d_j]$ if  $i\neq j$. 
The random graph  $\mathcal{G}_n(\alpha,\beta, W)$ is  therefore heterogeneous. 
When $w_{ij}$ $(1\leq i<j\leq n)$ are equal to one, $\mathcal{G}_n(\alpha,\beta, W)$ is the well-known Erd\"{o}s-R\'{e}nyi random graph. We simply denote it as $\mathcal{G}_n(\alpha)$.  The Erd\"{o}s-R\'{e}nyi random graph is homogeneous because the expected degrees of nodes are the same. The inhomogeneous Erd\"{o}s-R\'{e}nyi random graphs in \cite{CHHS20,Y23b,ZY23} are a special case  of $\mathcal{G}_n(\alpha,\beta, W)$. Moreover, the random graph $\mathcal{G}_n(\alpha,\beta, W)$ is studied in \cite{Y23,Y23c,YZ23}. We adopt $\mathcal{G}_n(\alpha,\beta, W)$ as the benchmark model in this paper.

The average closure coefficient of a graph is defined as
\begin{equation}\label{closurecoe}
\overline{\mathcal{H}}_n=\frac{1}{n}\sum_{i=1}^n\frac{\sum\limits_{j\neq k} A_{ij}A_{jk}A_{ki}}{\sum\limits_{j\neq k}A_{ij}A_{jk}},
\end{equation}
where any summation term with denominator zero is set to be zero \cite{YBL19}. The denominator in (\ref{closurecoe}) is the number of 2-path that starts at node $i$ and the numerator is 2 times the number of triangles involving node $i$. The difference between the average closure coefficient and the average clustering coefficient lies in the denominator.  The denominator in the average clustering coefficient is the number of wedges centered
at node $i$. In
real-world networks, the average closure coefficient and the average clustering coefficient can be positively correlated, negatively correlated, or weakly correlated. The local closure coefficient can therefore capture complementary information on fundamental clustering structure missed by the classical clustering coefficient \cite{YBL19}.

The theoretical properties of the average closure coefficient are not well studied. \cite{YBL19} obtained the expectation of the local average closure coefficient in the configuration model. In this paper, we derive the asymptotic distribution  of the average closure coefficient  $\overline{\mathcal{H}}_n$ in $\mathcal{G}_n(\alpha,\beta, W)$. This is not a straightforward task because $\overline{\mathcal{H}}_n$ is an average of dependent terms and each term is a ratio of dependent quantities.

\begin{Theorem}\label{mainthm} 
For the heterogeneous random graph $\mathcal{G}_n(\alpha,\beta, W)$, we have
    \[\frac{\overline{\mathcal{H}}_n-\mathbb{E}[\overline{\mathcal{H}}_n]}{\sigma_n}\Rightarrow \mathcal{N}(0,1),\]
where $\sigma_n^2=\sigma_{1n}^2+\sigma_{2n}^2$ and
\[\sigma_{1n}^2=\frac{4}{n^2}\sum_{i<j<k}\left(\frac{1}{\nu_i}+\frac{1}{\nu_j}+\frac{1}{\nu_k}\right)^2p_nw_{ij}(1-p_nw_{ij})p_nw_{jk}(1-p_nw_{jk})p_nw_{ki}(1-p_nw_{ki}),\]
\[\sigma_{2n}^2=\frac{1}{n^2}\sum_{i<j}\left(2b_{ij}+2c_{ij}+2c_{ji}-(a_i+a_j)-(e_{ij}+e_{ji})\right)^2p_nw_{ij}(1-p_nw_{ij}),\]
\[\nu_i=\sum_{j=1}^n\sum_{k=1}^np_n^2w_{ij}w_{jk}, \hskip 1cm b_{ij}=\sum_{k}\frac{p_nw_{ik}p_nw_{jk}}{\nu_k},\hskip 1cm c_{ij}=\sum_{k}\frac{p_nw_{ik}p_nw_{jk}}{\nu_i},\]
\[a_s=\sum_{i,j,k}\frac{p_nw_{ij}p_nw_{jk}p_nw_{ki}p_nw_{is}}{\nu_i^2},\hskip 1cm e_{is}=\sum_{j,k,t}\frac{p_nw_{ij}p_nw_{jk}p_nw_{ki}p_nw_{st}}{\nu_i^2}.\]

\end{Theorem}

\medskip

According to Theorem \ref{mainthm}, the standardized  average closure coefficient  converges in distribution to the standard normal distribution. Due to the dependence of the summation terms in $\overline{\mathcal{H}}_n$, the proof of Theorem \ref{mainthm} is not straightforward. Our proof strategy is as follows: expand each summation term to $k_0=\lceil \frac{4}{1-\alpha}\rceil +1$ order by Taylor expansion, single out the leading term and prove the remainder terms are negligible and the leading term converges in distribution to the standard normal distribution.

For
the Erd\"{o}s-R\'{e}nyi random graph $\mathcal{G}_n(\alpha)$, the variance $\sigma_n^2$ can be greatly simplified. We have the following corollary.

\begin{Corollary}\label{mcor}
\label{mainthmcor}
     For the Erd\"{o}s-R\'{e}nyi random graph $\mathcal{G}_n(\alpha)$, we have
    \[\frac{\overline{\mathcal{H}}_n-\mathbb{E}[\overline{\mathcal{H}}_n]}{\sigma_n}\Rightarrow \mathcal{N}(0,1),\]
    where $\sigma_n^2=\sigma_{1n}^2+\sigma_{2n}^2$, 
    $\sigma_{1n}^2=\frac{6}{n^{3-\alpha}}(1+o(1))$ and $\sigma_{2n}^2=\frac{2}{n^{2+\alpha}}(1+o(1))$. Hence we have
\begin{equation*} 
\sigma_n^2 =
    \begin{cases}
         \frac{6}{n^{3-\alpha}}(1+o(1)), & \text{if $\alpha>\frac{1}{2}$}, \\
   \frac{8}{n^2\sqrt{n}}(1+o(1)), & \text{if $\alpha=\frac{1}{2}$}, \\
    \frac{2}{n^{2+\alpha}}(1+o(1)) ,& \text{if $\alpha<\frac{1}{2}$}.     
    \end{cases}       
\end{equation*}
\end{Corollary}
For fixed large integer $n$, it is easy to verify that
\[\lim_{\alpha \to (\frac{1}{2})^{-}}\frac{2}{n^{2+\alpha}}=\frac{2}{n^{2}\sqrt{n}}\neq \frac{8}{n^2\sqrt{n}}, \]
\[\lim_{\alpha \to (\frac{1}{2})^{+}}\frac{6}{n^{3-\alpha}}=\frac{6}{n^{2}\sqrt{n}}\neq \frac{8}{n^2\sqrt{n}}. \]
As a function of $\alpha$, $\sigma_n^2$ does not change continuously  at $\alpha=\frac{1}{2}$. From this point of view, the scale $\sigma_n^2$ has a phase change at $\alpha=\frac{1}{2}$.

Interestingly, the leading order of $\sigma_n^2$ of the average closure coefficient in $\mathcal{G}_n(\alpha)$ is exactly the same as the average clustering coefficient \cite{YZ23}. The average closure coefficient and the average clustering coefficient exhibit the same phase change phenomenon.  In this sense, they are not significantly different in the Erd\"{o}s-R\'{e}nyi random graph $\mathcal{G}_n(\alpha)$.

\section{ Proof of main result}\label{main}

In this section, we provide detailed proof of Theorem \ref{mainthm}. Firstly, we present several lemmas. For convenience, denote $\mu_{ij}=p_nw_{ij}$, $V_i=\sum_{j,k\notin\{i\}}A_{ij}A_{jk}$ and $\bar{A}_{ij}=A_{ij}-\mu_{ij}$.

\begin{Lemma}\label{lem0}
 Let $\mathcal{G}_n(\alpha, \beta, W)$ be defined in Definition \ref{def1}, $\delta_n=\left(\log(np_n)\right)^{-2}$ and $M$ be a constant greater than $\frac{e^2}{1-p_n\beta}$. 
For any $i\in[n]$, we have
\[\mathbb{P}(d_{i}=k)\leq e^{-np_n\beta(1+o(1))},\ \ \ \ k\leq \delta_nnp_n,\]
\[\mathbb{P}(d_{i}=k)\leq e^{-np_n\beta(1+o(1))},\ \ \ \ k\geq Mnp_n.\]
\end{Lemma}

Lemma \ref{lem0} is proved in \cite{Y23b}. We omit the proof in this paper.

\begin{Lemma}\label{lem1}
Let $\epsilon_n=\left(\log(np_n)\right)^{-5}$ and $\delta_n=\left(\log(np_n)\right)^{-2}$. For the heterogeneous random graph $\mathcal{G}_n(\alpha,\beta, W)$, we have
    \[\mathbb{P}(V_1=k)\leq e^{-\beta np_n(1+o(1))}+e^{-\delta_n(np_n)^2\beta (1+o(1)},\ \ \ k\leq \epsilon_n (np_n)^2.\]
\end{Lemma}

\noindent
{\bf Proof of Lemma \ref{lem1}.} For each $j\neq 1$, denote $d_{j(1)}=\sum_{k\not\in\{j,1\}}A_{jk}$. Then $V_1=\sum_{j=1}^nA_{1j}d_{j(1)}$. Note that $A_{12}, A_{13}, \dots, A_{1n} $ are independent of $d_{2(1)}$, $d_{3(1)}$, $\dots$, $d_{n(1)}$. By the property of conditional probability and Lemma \ref{lem0}, we have
\begin{eqnarray}\nonumber
\mathbb{P}(V_1=k)&=&\sum_{t=1}^{n-1}\mathbb{P}(V_1=k|d_1=t)\mathbb{P}(d_1=t)\\ \nonumber
&=&\sum_{t=1}^{\delta_nnp_n}\mathbb{P}(V_1=k|d_1=t)\mathbb{P}(d_1=t)+\sum_{t=Mnp_n}^{n-1}\mathbb{P}(V_1=k|d_1=t)\mathbb{P}(d_1=t)\\ \nonumber
&&+\sum_{t=\delta_nnp_n+1}^{Mnp_n-1}\mathbb{P}(V_1=k|d_1=t)\mathbb{P}(d_1=t)\\ \label{lem1eq1}
&\leq&2e^{-\beta np_n(1+o(1))}+\sum_{t=\delta_nnp_n+1}^{Mnp_n-1}\mathbb{P}(V_1=k|d_1=t)\mathbb{P}(d_1=t).
\end{eqnarray}
Next we find an upper bound of the second term of (\ref{lem1eq1}).

Note that $d_1=t$ implies there are exactly $t$ of $A_{12}, A_{13}, \dots, A_{1n} $ are equal to 1 and $n-1-t$ of them are equal to zero. There are $\binom{n-1}{t}$ possible choices. Without loss of generality, let $A_{12}=A_{13}=\dots=A_{1(t+1)}=1$ and $A_{1(t+2)}=A_{1(t+3)}=\dots=A_{1n}=0$. Then 
\begin{eqnarray}\nonumber
&&\mathbb{P}(V_1=k|A_{12}=\dots=A_{1(t+1)}=1,A_{1(t+2)}=\dots=A_{1n}=0)\\ \nonumber
&=&\mathbb{P}(d_{2(1)}+d_{3(1)}+\dots+d_{(t+1)(1)}=k)\\ \label{lem1eq2}
&\leq&\mathbb{P}(d_{2(1)}+d_{3(1)}+\dots+d_{(t+1)(1)}\leq k).
\end{eqnarray}
It is easy to verify that
\begin{eqnarray*}
d_{2(1)}+d_{3(1)}+\dots+d_{(t+1)(1)}=2\sum_{2\leq i<j\leq t+1}A_{ij}+\sum_{i=2}^{t+1}\sum_{j=t+2}^nA_{ij}\geq\sum_{i=2}^{t+1}\sum_{j=t+2}^nA_{ij}.
\end{eqnarray*}
Then
\begin{eqnarray}\label{lem1eq3}
\mathbb{P}(d_{2(1)}+d_{3(1)}+\dots+d_{(t+1)(1)}\leq k)\leq\mathbb{P}\left(\sum_{i=2}^{t+1}\sum_{j=t+2}^nA_{ij}\leq k\right).
\end{eqnarray}

Let $N_t=\{(i,j)|2\leq i\leq t+1,  t+2\leq j\leq n\}$ and $\theta_{t}=\{p_nw_{ij}|(i,j)\in N_t\}$. Then
$\sum_{i=2}^{t+1}\sum_{j=t+2}^nA_{ij}$ follows the Poisson-Binomial distribution $PB(\theta_{t})$. Recall that $\beta\leq w_{ij}\leq1$. Then 
\begin{eqnarray}\nonumber
\mathbb{P}\left(\sum_{i=2}^{t+1}\sum_{j=t+2}^nA_{ij}=k\right)&=&\sum_{S\subset N_t,|S|=k}\prod_{(i,j)\in S}p_nw_{ij}\prod_{(i,j)\in S^C}(1-p_nw_{ij})\\ \nonumber
&\leq&\sum_{S\subset N_t,|S|=k}\prod_{(i,j)\in S}p_n\prod_{(i,j)\in S^C}(1-p_n\beta)\\ \label{lem1eq4}
&\leq &\binom{nt}{k}p_n^k(1-p_n\beta)^{nt-t^2-t-k}.
\end{eqnarray}
Note that $\binom{nt}{k}\leq e^{k\log (nt)-k\log k+k}$ and 
$(1-p_n\beta)^{nt-t^2-t-k}=e^{(nt-t^2-t-k)\log(1-p_n\beta)}$.
Then by (\ref{lem1eq4}) we get
\begin{eqnarray}\nonumber
\mathbb{P}\left(\sum_{i=2}^{t+1}\sum_{j=t+2}^nA_{ij}=k\right)
\leq e^{g(k)},
\end{eqnarray}
where $g(k)=k\log (ntp_n)-k\log k+k+(nt-t^2-t-k)\log(1-p_n\beta)$. Considering $k$ as continuous variable, the derivative of $g(k)$ with respect to $k$ is equal to 
\[g^{\prime}(k)=\log\left(\frac{ntp_n}{1-p_n\beta}\right)-\log k.\]
Clearly, $g^{\prime}(k)>0$ for $k<\frac{ntp_n}{1-p_n\beta}$ and $g^{\prime}(k)<0$ for $k>\frac{ntp_n}{1-p_n\beta}$. Then $g(k)$ achieves its maximum at $k=\frac{ntp_n}{1-p_n\beta}$. Let $c_n=\left(\log(np_n)\right)^{-\frac{1}{2}}$. For $k\leq c_nntp_n$,  $g(k)\leq g(c_nntp_n)$. Note that $-c_n\log(c_n)=o(1)$. Hence
\begin{eqnarray*}\nonumber
\mathbb{P}\left(\sum_{i=2}^{t+1}\sum_{j=t+2}^nA_{ij}=k\right)
&\leq&e^{c_nntp_n\log\frac{1}{c_n(1-p_n\beta)}+c_nntp_n+nt\log(1-p_n\beta)}e^{-\frac{t^2+t}{nt}\log(1-p_n\beta)}\\ 
&\leq& e^{-ntp_n\beta(1+o(1))}.
\end{eqnarray*}

Note that $\delta_n np_n\leq t\leq Mnp_n$ in the second term of (\ref{lem1eq1}).  Then  
\[k\leq \epsilon_n(np_n)^2= \frac{(np_n)^2}{\left(\log (np_n)\right)^5}\leq  \frac{(np_n)^2}{\left(\log (np_n)\right)^2\sqrt{\log(np_n)}} \leq c_nntp_n.\]
Hence, for $k\leq \epsilon_n(np_n)^2$, we have
\begin{eqnarray}\label{lem1eq5}
\mathbb{P}\left(\sum_{i=2}^{t+1}\sum_{j=t+2}^nA_{ij}=k\right)
\leq e^{-\delta_n (np_n)^2\beta(1+o(1))}.
\end{eqnarray}

Combining (\ref{lem1eq1})- (\ref{lem1eq5}), one has
\begin{eqnarray}\nonumber
\sum_{t=\delta_nnp_n+1}^{Mnp_n-1}\mathbb{P}(V_1=k|d_1=t)\mathbb{P}(d_1=t)&\leq& Mnp_n \binom{n-1}{t} ke^{-\delta_n (np_n)^2\beta(1+o(1))}\\ \label{lem1eq6}
&=&e^{-\delta_n (np_n)^2\beta(1+o(1))}.
\end{eqnarray}
Based on  (\ref{lem1eq1}) and (\ref{lem1eq6}), the result of Lemma \ref{lem1} holds.



\qed

\begin{Lemma}\label{lem2}
For the heterogeneous random graph $\mathcal{G}_n(\alpha,\beta, W)$ and a positive integer $t$, we have
    \[\mathbb{E}[(V_i-\nu_i)^{2t}]=O\left((np_n)^{3t}\right),\]
    uniformly for all $i\in[n]$.
\end{Lemma}

\noindent
{\bf Proof of Lemma \ref{lem2}.} It is straightforward to get that
\begin{eqnarray}\label{lem2eq1}
    V_i-\nu_i=\sum_{j\neq k}\bar{A}_{ij}\bar{A}_{jk}+\sum_{j\neq k}\bar{A}_{ij}\mu_{jk}+\sum_{j\neq k}\mu_{ij}\bar{A}_{jk}.
\end{eqnarray}
Then
\begin{eqnarray}\label{lem2eq2}
    (V_i-\nu_i)^{2t}\leq 3^{2t}\left[\left(\sum_{j\neq k}\bar{A}_{ij}\bar{A}_{jk}\right)^{2t}+\left(\sum_{j\neq k}\bar{A}_{ij}\mu_{jk}\right)^{2t}+\left(\sum_{j\neq k}\mu_{ij}\bar{A}_{jk}\right)^{2t}\right].
\end{eqnarray}
Next we find upper bound of the expectation of each term in (\ref{lem2eq2}).

Consider the first term in (\ref{lem2eq2}).
Note that
\begin{eqnarray}\label{lem2eq3}
    \mathbb{E}\left[\left(\sum_{j\neq k}\bar{A}_{ij}\bar{A}_{jk}\right)^{2t}\right]
    &=&\sum_{\substack{j_1\neq k_1, \dots,j_{2t}\neq k_{2t}} }\mathbb{E}\left[\bar{A}_{ij_1}\bar{A}_{j_1k_1}\dots \bar{A}_{ij_{2t}}\bar{A}_{j_{2t}k_{2t}}\right].
\end{eqnarray}

Note that $\mathbb{E}[\bar{A}_{ij}\bar{A}_{st}]=0$ if $\{i,j\}\neq \{s,t\}$ and $\mathbb{E}[|\bar{A}_{ij}^t|]=O(p_n)$ for any positive integer $t$. If $\{i,j_1\}\neq\{i,j_{s}\}$ for all $s\geq2$ and $\{i,j_1\}\neq \{j_{l},k_{l}\}$ for all $l\in\{1,2,\dots,2t\}$, then
\[\mathbb{E}\left[\bar{A}_{ij_1}\bar{A}_{j_1k_1}\dots \bar{A}_{ij_{2t}}\bar{A}_{j_{2t}k_{2t}}\right]=\mathbb{E}\left[\bar{A}_{ij_1}\right]\mathbb{E}\left[\bar{A}_{j_1k_1}\dots \bar{A}_{ij_{2t}}\bar{A}_{j_{2t}k_{2t}}\right]=0.\]
Hence, for the expectation in (\ref{lem2eq3}) to be non-zero, $\{i,j_1\}=\{i,j_{s}\}$ for some $s\geq2$ or $\{i,j_1\}=\{j_{l},k_{l}\}$ for some $l\in\{1,2,\dots,2t\}$. Note that $j_{l}\neq i$ and $k_{l}\neq i$ for all $l$. It is impossible that $\{i,j_1\}=\{j_{l},k_{l}\}$ for some $l\in\{1,2,\dots,2t\}$. Hence $\{i,j_1\}=\{i,j_{s}\}$ for some $s\geq2$. Similarly, for each $s\in\{1,2,\dots,2t\}$, there exists $s_1\in\{1,2,\dots,s-1,s+1,\dots,2t\}$ such that $\{i,j_s\}=\{i,j_{s_1}\}$. Hence, $|\{j_1,j_2,\dots,j_{2t}\}|\leq t$. Without loss of generality, assume   $\{j_1,j_2,\dots,j_{2t}\}=\{j_1,j_2,\dots,j_s\}$ for some $s\leq t$ and $|\{j_1,j_2,\dots,j_s\}|=s$. Let $t_l$ be the number of elements in $\{j_1,j_2,\dots,j_{2t}\}$ that are equal to $j_l$, $1\leq l\leq s$. Then $t_1+t_2+\dots+t_s=2t$. Suppose
$j_{r_{lq}}=j_l$ for $q=1,2,\dots,t_l$, $1\leq l\leq s$ and $r_{lq}\in \{1,2,\dots,2t\}$. In this case,
\begin{eqnarray}\nonumber
\mathbb{E}\left[\bar{A}_{ij_1}\bar{A}_{j_1k_1}\dots \bar{A}_{ij_{2t}}\bar{A}_{j_{2t}k_{2t}}\right]&=& \mathbb{E}\left[\left(\prod_{l=1}^s\mathbb{E}\left[\bar{A}_{ij_{l}}^{t_l}\right]\right)\left(\prod_{l=1}^s\prod_{q=1}^{t_{l}}\bar{A}_{j_{l}k_{r_{lq}}}\right]\right)\\ \label{lem2eq4}
&=&O(p_n^s)\mathbb{E}\left[\prod_{l=1}^s\prod_{q=1}^{t_{l}}\bar{A}_{j_{l}k_{r_{lq}}}\right].
\end{eqnarray}

If $\bar{A}_{j_{1}k_{r_{11}}}\neq \bar{A}_{j_{l_1}k_{r_{l_1q_1}}}$ for any $l_1,q_1$ with $l_1\neq 1$ or $q_1\neq1$, then  
\[\mathbb{E}\left[\prod_{l=1}^s\prod_{q=1}^{t_{l}}\bar{A}_{j_{l}k_{r_{lq}}}\right]=\mathbb{E}[\bar{A}_{j_{1}k_{r_{11}}}]\mathbb{E}\left[\left(\prod_{q=2}^{t_{1}}\bar{A}_{j_{1}k_{r_{1q}}}\right)\prod_{l=2}^s\prod_{q=1}^{t_{l}}\bar{A}_{j_{l}k_{r_{lq}}}\right]=0.\]
 Hence $\bar{A}_{j_{1}k_{r_{11}}}=\bar{A}_{j_{l_1}k_{r_{l_1q_1}}}$ for some $l_1,q_1$ with $l_1\neq 1$ or $q_1\neq1$. Similarly, for each $l,q$, there exist $l_1,q_1$ with $l_1\neq l$ or $q_1\neq q$  such that  $\bar{A}_{j_{l}k_{r_{lq}}}=\bar{A}_{j_{l_1}k_{r_{l_1q_1}}}$. That is, $(j_{l},k_{r_{lq}})=(j_{l_1},k_{r_{l_1q_1}})$. Otherwise, the expectation in (\ref{lem2eq4}) is zero. In this case, one has
 \begin{equation}\label{lem2eq5}
    \Big| \Big\{(j_1,k_{r_{11}}),(j_1,k_{r_{12}}),\dots,(j_1,k_{r_{1t_1}}),\dots, (j_s,k_{r_{s1}}),\dots, (j_s,k_{r_{st_s}})\Big\}\Big|=m\leq t.
 \end{equation}
 Then
 \begin{equation}\label{lem2eq05}
\mathbb{E}\left[\prod_{l=1}^s\prod_{q=1}^{t_{l}}\bar{A}_{j_{l}k_{r_{lq}}}\right]=O(p_n^m).
 \end{equation}
Note that $(j_{l},k_{r_{lq}})=(j_{l_1},k_{r_{l_1q_1}})$ with $l_1\neq l$ or $q_1\neq q$ implies $|\{j_{l},k_{r_{lq}},j_{l_1},k_{r_{l_1q_1}}\}|\leq 3$. That is, if one pair $(j_{l},k_{r_{lq}})$ is equal to another pair, the number of distinct indices will reduce by at least 1. By (\ref{lem2eq5}), the number of pairs $(j_{l},k_{r_{lq}})$ is reduced by $2t-m$. Hence, the number of distinct indices $j_{l},k_{r_{lq}}$ is reduced by at least $2t-m$. Then
 \begin{equation}\label{lem2eq6}
    \Big| \Big\{j_1,j_2,\dots,j_s,k_{r_{11}},k_{r_{12}},\dots,k_{r_{1t_1}},\dots, k_{r_{s1}},\dots, k_{r_{st_s}}\Big\}\Big|\leq (s+2t)-(2t-m)=s+m.
 \end{equation}
There are at most $n^{s+m}$ choices of indices $j_l,k_r$ satisfying (\ref{lem2eq6}). Combining (\ref{lem2eq3}), (\ref{lem2eq4}) and (\ref{lem2eq05}) yields
\begin{eqnarray}\label{lem2eq7}
    \mathbb{E}\left[\left(\sum_{j\neq k}\bar{A}_{ij}\bar{A}_{jk}\right)^{2t}\right]=O\left((np_n)^{s+m}\right)=O\left((np_n)^{2t}\right).
\end{eqnarray}
By a similar proof of Lemma 3.2 in \cite{YZ23}, it is easy to get 
\begin{eqnarray}\label{lem2eq8}
    \mathbb{E}\left[\left(\sum_{j\neq k}\bar{A}_{ij}\mu_{jk}\right)^{2t}\right]=O\left((np_n)^{3t}\right),
\end{eqnarray}
\begin{eqnarray}\label{lem2eq9}
    \mathbb{E}\left[\left(\sum_{j\neq k}\mu_{ij}\bar{A}_{jk}\right)^{2t}\right]=O\left((np_n)^{2t}p_n^t\right).
\end{eqnarray}
By (\ref{lem2eq2}), (\ref{lem2eq7}), (\ref{lem2eq8}) and (\ref{lem2eq9}), the proof is complete.

\qed

\subsection{Proof of Theorem \ref{mainthm}}
For convenience, denote $\Delta_{ijk}=A_{ij}A_{jk}A_{ki}$ and  $\Delta_{i}=\sum_{j,k}A_{ij}A_{jk}A_{ki}$. Let $k_0=\lceil \frac{4}{1-\alpha}\rceil +1$. Applying Taylor expansion to the function $f(x)=\frac{1}{x}$ at $\nu_i=\sum_{j,k}\mu_{ij}\mu_{jk}$ yields
\begin{eqnarray*}\nonumber
\overline{\mathcal{H}}&=&\frac{1}{n}\sum_{i=1}^n\frac{\Delta_{i}}{V_i}\\ 
&=&\frac{1}{n}\sum_{i=1}^n\frac{\Delta_{i}}{\nu_i}-\frac{1}{n}\sum_{i=1}^n\frac{\Delta_{i}(V_i-\nu_i)}{\nu_i^2}+\sum_{t=2}^{k_0-1}(-1)^t\frac{1}{n}\sum_{i=1}^n\frac{\Delta_{i}(V_i-\nu_i)^t}{\nu_i^{t+1}}+(-1)^{k_0}\frac{1}{n}\sum_{i=1}^n\frac{\Delta_{i}(V_i-\nu_i)^{k_0}}{X_i^{k_0+1}}.
\end{eqnarray*}
where $X_i$ is between $V_i$ and $\nu_i$. Then we have
\begin{eqnarray}\nonumber
\overline{\mathcal{H}}-\mathbb{E}\left[\overline{\mathcal{H}}\right]&=&\frac{1}{n}\sum_{i=1}^n\frac{\Delta_{i}-\mathbb{E}\left[\Delta_{i}\right]}{\nu_i}-\frac{1}{n}\sum_{i=1}^n\frac{\Delta_{i}(V_i-\nu_i)-\mathbb{E}\left[\Delta_{i}(V_i-\nu_i)\right]}{\nu_i^2}\\ \nonumber
&&+\sum_{t=2}^{k_0-1}(-1)^{t}\frac{1}{n}\sum_{i=1}^n\frac{\Delta_{i}(V_i-\nu_i)^t-\mathbb{E}\left[\Delta_{i}(V_i-\nu_i)^t\right]}{\nu_i^{t+1}}\\ \label{eq1}
&&+(-1)^{k_0}\frac{1}{n}\sum_{i=1}^n\frac{\Delta_{i}(V_i-\nu_i)^{k_0}}{X_i^{k_0+1}}-(-1)^{k_0}\frac{1}{n}\sum_{i=1}^n\mathbb{E}\left[\frac{\Delta_{i}(V_i-\nu_i)^{k_0}}{X_i^{k_0+1}}\right].
\end{eqnarray}
Next we prove the first two terms of (\ref{eq1}) are the leading terms and the last three terms are negligible.

\noindent
{\bf Consider the first term of (\ref{eq1})}. By the proof of \cite{YZ23}, we have the following result. If $\alpha>\frac{1}{2}$, then
\begin{eqnarray}\label{y1ng}
\frac{1}{n}\sum_{i=1}^n\frac{\Delta_{i}-\mathbb{E}\left[\Delta_{i}\right]}{\nu_i}
&=&\frac{1}{n}\sum_{i\neq j\neq k}\frac{\bar{A}_{ij}\bar{A}_{ik}\bar{A}_{jk}}{\nu_i}+o_P\left(\frac{1}{n\sqrt{np_n}}\right).
\end{eqnarray}
If $\alpha<\frac{1}{2}$, then
\begin{eqnarray} \label{y1ns}
\frac{1}{n}\sum_{i=1}^n\frac{\Delta_{i}-\mathbb{E}\left[\Delta_{i}\right]}{\nu_i}
=\frac{1}{n}\sum_{i\neq j\neq k}\frac{\left(\bar{A}_{jk}\mu_{ij}\mu_{ik}+\bar{A}_{ij}\mu_{jk}\mu_{ik}+\bar{A}_{ik}\mu_{jk}\mu_{ij}\right)}{\nu_i}+o_P\left(\frac{\sqrt{p_n}}{n}\right).
\end{eqnarray}
If $\alpha=\frac{1}{2}$, then
\begin{eqnarray} \nonumber
&&\frac{1}{n}\sum_{i=1}^n\frac{\Delta_{i}-\mathbb{E}\left[\Delta_{i}\right]}{\nu_i}\\ \label{y1eq}
&=&\frac{1}{n}\sum_{i\neq j\neq k}\frac{\left(\bar{A}_{ij}\bar{A}_{ik}\bar{A}_{jk}+\bar{A}_{jk}\mu_{ij}\mu_{ik}+\bar{A}_{ij}\mu_{jk}\mu_{ik}+\bar{A}_{ik}\mu_{jk}\mu_{ij}\right)}{\nu_i}+o_P\left(\frac{\sqrt{p_n}}{n}\right).
\end{eqnarray}

\medskip

\noindent
{\bf Consider the last two terms of (\ref{eq1})}.  Let $\epsilon_n$ be defined in Lemma \ref{lem2}. It is easy to get that
\begin{eqnarray}\nonumber
\frac{1}{n}\sum_{i=1}^n\mathbb{E}\left[\left|\frac{\Delta_{i}(V_i-\nu_i)^{k_0}}{X_i^{k_0+1}}\right|\right]&=&\frac{1}{n}\sum_{i=1}^n\mathbb{E}\left[\left|\frac{\Delta_{i}(V_i-\nu_i)^{k_0}}{X_i^{k_0+1}}\right|I[X_i\geq \epsilon_n (np_n)^2]\right]\\ \label{eqremin}
&&+\frac{1}{n}\sum_{i=1}^n\mathbb{E}\left[\left|\frac{\Delta_{i}(V_i-\nu_i)^{k_0}}{X_i^{k_0+1}}\right|I[X_i<\epsilon_n (np_n)^2]\right].
\end{eqnarray}

Note that $k_0=\lceil \frac{4}{1-\alpha}\rceil +1>\frac{4}{1-\alpha}$.
By Lemma \ref{lem2} and the Cauchy-Schwarz inequality, we have
\begin{eqnarray}\nonumber
&&\frac{1}{n}\sum_{i=1}^n\mathbb{E}\left[\left|\frac{\Delta_{i}(V_i-\nu_i)^{k_0}}{X_i^{k_0+1}}\right|I[X_i\geq \epsilon_n (np_n)^2]\right]\\ \nonumber
&\leq&\frac{1}{n\epsilon_n^{k_0+1} (np_n)^{2k_0+2}}\sum_{i=1}^n\sqrt{\mathbb{E}[\Delta_{i}^2]\mathbb{E}[(V_i-\nu_i)^{2k_0}]}\\ \nonumber
&=&O\left(\frac{p_n(np_n)^{\frac{3}{2}k_0+2}}{\epsilon_n^{k_0+1} (np_n)^{2k_0+2}}\right)\\ \label{eqgreat}
&=&o\left(\frac{p_n}{n^2}\right).
\end{eqnarray}

Recall that $\nu_i=\Theta((np_n)^2)$ and $X_i$ is between $\nu_i$ and $V_i$. When $X_i<\epsilon_n (np_n)^2$, we have $V_i<X_i$. Otherwise, $X_i$ can not be between  $\nu_i$ and $V_i$. Moreover, if $V_i=0$, then the $i$-th term in the definition of the average clustering coefficient in (\ref{closurecoe}) vanishes. Hence, we assume $V_i\geq1$ and then $X_i\geq1$. It is easy to verify that $\Delta_i\leq n^2$ and $|V_i-\nu_i|\leq n^2$. By Lemma \ref{lem1}, we have
\begin{eqnarray}\nonumber
&&\frac{1}{n}\sum_{i=1}^n\mathbb{E}\left[\left|\frac{\Delta_{i}(V_i-\nu_i)^{k_0}}{X_i^{k_0+1}}\right|I[X_i<\epsilon_n (np_n)^2]\right]\\ \nonumber
&=&O\left(n^{2k_0+2}\right)\max_{1\leq i\leq n}\mathbb{P}(V_i<\epsilon_n (np_n)^2)\\ \nonumber
&\leq&n^{2k_0+2}\sum_{k=1}^{\epsilon_n (np_n)^2}\mathbb{P}(V_1=k)\\ \label{eqsmall}
&=&e^{-np_n\beta (1+o(1)}.
\end{eqnarray}

Combining (\ref{eqremin}), (\ref{eqgreat}) and (\ref{eqsmall}) yields
\begin{eqnarray}\label{eqresmall}
\frac{1}{n}\sum_{i=1}^n\frac{\Delta_{i}(V_i-\nu_i)^{k_0}}{X_i^{k_0+1}}-\frac{1}{n}\sum_{i=1}^n\mathbb{E}\left[\frac{\Delta_{i}(V_i-\nu_i)^{k_0}}{X_i^{k_0+1}}\right]=o_P\left(\frac{1}{n\sqrt{np_n}}+\frac{\sqrt{p_n}}{n}\right).
\end{eqnarray}

\medskip

\noindent
{\bf Consider the second term of (\ref{eq1})}. By (\ref{lem2eq1}), it is easy to verify that
\begin{eqnarray}\nonumber
\frac{1}{n}\sum_{i=1}^n\frac{\Delta_{i}(V_i-\nu_i)}{\nu_i^2}&=&\frac{1}{n}\sum_{i\neq j\neq k, s\neq t}\frac{A_{ij}A_{jk}A_{ki}\bar{A}_{is}\bar{A}_{st}}{\nu_i^2}+\frac{1}{n}\sum_{i\neq j\neq k, s\neq t}\frac{A_{ij}A_{jk}A_{ki}\bar{A}_{is}\mu_{st}}{\nu_i^2}\\ \label{seceq1}
&&+\frac{1}{n}\sum_{i\neq j\neq k, s\neq t}\frac{A_{ij}A_{jk}A_{ki}\mu_{is}\bar{A}_{st}}{\nu_i^2}.
\end{eqnarray}
Next we show the first term of (\ref{seceq1}) is negligible and the last two terms are leading terms. 

Consider the first term of (\ref{seceq1}).
If $\{j,k\}=\{s,t\}$ in the summation of the first term of (\ref{seceq1}), then
\begin{eqnarray}\label{seteq1}
\mathbb{E}\left[\frac{1}{n}\sum_{i\neq j\neq k}\left|\frac{A_{ij}A_{jk}A_{ki}\bar{A}_{ij}\bar{A}_{jk}}{\nu_i^2}\right|\right]=O\left(\frac{1}{n(np_n)}\right).
\end{eqnarray}
Suppose $|\{j,k\}\cap\{s,t\}|=1$ in the summation of the first term of (\ref{seceq1}). Without loss of generality, let $s=j$ and $t\neq k$. Then
\begin{eqnarray}\nonumber
&&\mathbb{E}\left[\left(\frac{1}{n}\sum_{i\neq j\neq k\neq t}\frac{A_{ij}A_{jk}A_{ki}\bar{A}_{ij}\bar{A}_{jt}}{\nu_i^2}\right)^2\right]\\ \label{seteq2}
&=&\frac{1}{n^2}\sum_{\substack{i\neq j\neq k\neq t\\  i_1\neq j_1\neq k_1\neq t_1}}\frac{\mathbb{E}\left[A_{ij}A_{jk}A_{ki}\bar{A}_{ij}\bar{A}_{jt}A_{i_1j_1}A_{j_1k_1}A_{k_1i_1}\bar{A}_{i_1j_1}\bar{A}_{j_1t_1}\right]}{\nu_i^2\nu_{i_1}^2} .
\end{eqnarray}
If $(j,t)\notin\{(i_1,j_1),(j_1,k_1),(k_1,i_1),(j_1,t_1)\}$, then
 $\bar{A}_{jt}$ is independent of $A_{i_1j_1}$, $A_{i_1j_1}$, $A_{j_1k_1}$, $A_{k_1i_1}$. In this case,
\begin{eqnarray*}
&&\mathbb{E}\left[A_{ij}A_{jk}A_{ki}\bar{A}_{ij}\bar{A}_{jt}A_{i_1j_1}A_{j_1k_1}A_{k_1i_1}\bar{A}_{i_1j_1}\bar{A}_{j_1t_1}\right]\\
&=&\mathbb{E}[\bar{A}_{jt}]\mathbb{E}\left[A_{ij}A_{jk}A_{ki}\bar{A}_{ij}A_{i_1j_1}A_{j_1k_1}A_{k_1i_1}\bar{A}_{i_1j_1}\bar{A}_{j_1t_1}\right]=0.
\end{eqnarray*}
 Hence $(j,t)\in\{(i_1,j_1),(j_1,k_1),(k_1,i_1),(j_1,t_1)\}$. Similarly, $(j_1,t_1)\in\{(i,j),(j,k),(k,i),(j,t)\}$. Denote $e_1=(i,j)$, $e_2=(j,k)$, $e_3=(k,i)$, $e_4=(j,t)$,$e_5=(i_1,j_1)$, $e_6=(j_1,k_1)$, $e_7=(k_1,i_1)$, $e_8=(j_1,t_1)$,
 \[E=\Big\{e_1,e_2,e_3,e_4,e_5,e_6,e_7,e_8\Big\},\hskip 1cm F=\{i,j,k,t,i_1,j_1,k_1,t_1\}.\]
Then $|E|\leq 7$. If $|E|=7$, then $(j_1,t_1)=(j,t)$, $|F|=6$.
There are at most $n^{6}$ choices of the indices in $F$. Suppose $|E|\leq 6$. Let $e_{l_1}\in E $  and $e_{l_2}\in E $ with $l_1\neq l_2$. If $e_{l_1}=e_{l_2}$, then $|F|$ will reduce by at least 1. Hence, $|F|\leq |E|\leq 6$. In this case, there are at most $n^{|E|}$ choices of the indices in $F$. Note that
\[\mathbb{E}\left[A_{ij}A_{jk}A_{ki}\bar{A}_{ij}\bar{A}_{jt}A_{i_1j_1}A_{j_1k_1}A_{k_1i_1}\bar{A}_{i_1j_1}\bar{A}_{j_1t_1}\right]=O(p_n^{|E|}).\]
By (\ref{seteq2}) we conclude
\begin{eqnarray}\label{seteq3}
\mathbb{E}\left[\left(\frac{1}{n}\sum_{i\neq j\neq k\neq t}\frac{A_{ij}A_{jk}A_{ki}\bar{A}_{ij}\bar{A}_{jt}}{\nu_i^2}\right)^2\right]
=O\left(\frac{(np_n)^6}{n^2(np_n)^8}\right)=O\left(\frac{1}{n^2(np_n)^2}\right).
\end{eqnarray}

Suppose $|\{j,k\}\cap\{s,t\}|=0$ in the summation of the first term of (\ref{seceq1}). Then
\begin{eqnarray} \nonumber
&&\mathbb{E}\left[\left(\frac{1}{n}\sum_{i\neq j\neq k\neq t\neq s}\frac{A_{ij}A_{jk}A_{ki}\bar{A}_{is}\bar{A}_{st}}{\nu_i^2}\right)^2\right]\\  \label{0seteq4}
&=&\frac{1}{n^2}\sum_{\substack{i\neq j\neq k\neq t\neq s\\ i_1\neq j_1\neq k_1\neq t_1\neq s_1}}\frac{\mathbb{E}\left[A_{ij}A_{jk}A_{ki}\bar{A}_{is}\bar{A}_{st}A_{i_1j_1}A_{j_1k_1}A_{k_1i_1}\bar{A}_{i_1s_1}\bar{A}_{s_1t_1}\right]}{\nu_i^2\nu_{i_1}^2}.
\end{eqnarray}
By a similar argument as in (\ref{seteq3}), $(i,s),(s,t)\in\{(i_1,j_1),(j_1,k_1),(k_1,i_1),(i_1,s_1),(s_1,t_1)\}$ and $(i_1,j_1),(j_1,t_1)\in\{(i,j),(j,k),(k,i),(i,s),(s,t)\}$. Denote
 \[E=\Big\{(i,j),(j,k),(k,i),(i,s),(s,t),(i_1,j_1),(j_1,k_1),(k_1,i_1),(i_1,s_1),(s_1,t_1)\Big\}.\]
 \[F=\{i,j,k,s,t,i_1,j_1,k_1,s_1,t_1\}\]
Then $|E|\leq 8$. If $|E|=8$, then $\{(i_1,s_1),(s_1,t_1)\}=\{(i,j),(j,t)\}$ and $|F|=7$. There are at most $n^{7}$ choices of the indices in $F$. Suppose $|E|=7$. Then  $\{(i_1,s_1),(s_1,t_1),e_1\}=\{(i,j),(j,t),e\}$ for some  $e_1\in\{(i_1,j_1),(j_1,k_1),(k_1,i_1)\}$  and some $e\in\{(i,j),(j,k),(k,i)\}$. Then $|F|\leq 6$. Suppose $|E|\leq 6$. Then $|F|\leq |E|\leq 6$. In this case, there are at most $n^{|E|}$ choices of the indices in $F$. Note that
\[
\mathbb{E}\left[A_{ij}A_{jk}A_{ki}\bar{A}_{is}\bar{A}_{st}A_{i_1j_1}A_{j_1k_1}A_{k_1i_1}\bar{A}_{i_1s_1}\bar{A}_{s_1t_1}\right]=O(p_n^{|E|}).\]
By (\ref{0seteq4}) we conclude
\begin{eqnarray}\label{seteq4}
\mathbb{E}\left[\left(\frac{1}{n}\sum_{i\neq j\neq k\neq t\neq s}\frac{A_{ij}A_{jk}A_{ki}\bar{A}_{is}\bar{A}_{st}}{\nu_i^2}\right)^2\right]=O\left(\frac{n^7p_n^8+n^6p_n^6}{n^2(np_n)^8}\right)=o\left(\frac{1}{n^2(np_n)}\right).
\end{eqnarray}

Combining (\ref{seteq1})-(\ref{seteq4}) yields
\begin{eqnarray} \label{seteq5}
\frac{1}{n}\sum_{i\neq j\neq k, s\neq t}\frac{A_{ij}A_{jk}A_{ki}\bar{A}_{is}\bar{A}_{st}}{\nu_i^2}-
\mathbb{E}\left[\frac{1}{n}\sum_{i\neq j\neq k, s\neq t}\frac{A_{ij}A_{jk}A_{ki}\bar{A}_{is}\bar{A}_{st}}{\nu_i^2}\right]=o_P\left(\frac{1}{n\sqrt{np_n}}\right).
\end{eqnarray}

Consider the second term of (\ref{seceq1}). The second term of (\ref{seceq1}) can be expressed as
\begin{eqnarray}\nonumber
\frac{1}{n}\sum_{i\neq j\neq k, s\neq t}\frac{A_{ij}A_{jk}A_{ki}\bar{A}_{is}\mu_{st}}{\nu_i^2}
&=&\frac{1}{n}\sum_{i\neq j\neq k\neq s,t}\frac{A_{ij}A_{jk}A_{ki}\bar{A}_{is}\mu_{st}}{\nu_i^2}+\frac{1}{n}\sum_{i\neq j\neq k,t}\frac{A_{ij}\bar{A}_{ij}A_{jk}A_{ki}\mu_{jt}}{\nu_i^2}\\ \label{seteq6}
&&+\frac{1}{n}\sum_{i\neq j\neq k,t}\frac{A_{ij}A_{jk}A_{ki}\bar{A}_{ik}\mu_{kt}}{\nu_i^2}.
\end{eqnarray}
Now we find the order of each term in (\ref{seteq6}).

By a similar argument as in (\ref{seteq3}), we have
\begin{eqnarray}\label{seteq08}
\mathbb{E}\left[\left(\frac{1}{n}\sum_{i\neq j\neq k\neq s,t}\frac{A_{ij}A_{jk}A_{ki}\bar{A}_{is}\mu_{st}}{\nu_i^2}\right)^2\right]=O\left(\frac{p_n}{n^2}\right).
\end{eqnarray}

The variance of the second term of (\ref{seteq6}) is equal to
\begin{eqnarray}\nonumber
&&\mathbb{E}\left[\left(\frac{1}{n}\sum_{i\neq j\neq k,t}\frac{A_{ij}\bar{A}_{ij}A_{jk}A_{ki}\mu_{jt}-\mathbb{E}\left[A_{ij}\bar{A}_{ij}A_{jk}A_{ki}\mu_{jt}\right]}{\nu_i^2}\right)^2\right]\\ \nonumber
&=&\frac{1}{n^2}\sum_{\substack{i\neq j\neq k,t\\ i_1\neq j_1\neq k_1,t_1}}\frac{1}{\nu_i^2\nu_{i_1}^2}\mathbb{E}\Bigg[(A_{ij}\bar{A}_{ij}A_{jk}A_{ki}\mu_{jt}-\mathbb{E}\left[A_{ij}\bar{A}_{ij}A_{jk}A_{ki}\mu_{jt}\right])\\ \label{seteq8}
&&\times (A_{i_1j_1}\bar{A}_{i_1j_1}A_{j_1k_1}A_{k_1i_1}\mu_{j_1t_1}-\mathbb{E}\left[A_{i_1j_1}\bar{A}_{i_1j_1}A_{j_1k_1}A_{k_1i_1}\mu_{j_1t_1}\right])\Bigg].
\end{eqnarray}
Let $S=\{(i,j),(j,k),(k,i)\}\cap\{(i_1,j_1),(j_1,k_1),(k_1,i_1)\}$.
If $S=\emptyset$, then the expectation in (\ref{seteq8}) is zero. If $|S|=1$, then $|\{i,j,k,i_1,j_1,k_1\}|=4$. Denote $B_{ijkt}=A_{ij}\bar{A}_{ij}A_{jk}A_{ki}\mu_{jt}-\mathbb{E}\left[A_{ij}\bar{A}_{ij}A_{jk}A_{ki}\mu_{jt}\right]$. Then
\[\mathbb{E}\Big[B_{ijkt}B_{i_1j_1k_1t_1}\Big]=O(p_n^7).\]
Hence we have
\begin{eqnarray}\label{seteq9}
\frac{1}{n^2}\sum_{\substack{i\neq j\neq k,t\\ i_1\neq j_1\neq k_1,t_1\\ |S|=1}}\frac{1}{\nu_i^2\nu_{i_1}^2}\mathbb{E}\Big[B_{ijkt}B_{i_1j_1k_1t_1}\Big]
=O\left(\frac{n^6p_n^7}{n^2(np_n)^8}\right) 
=O\left(\frac{1}{n^2(np_n)^2}\right).
\end{eqnarray}
If
$|S|=2$, then $|\{i,j,k,i_1,j_1,k_1\}|=3$  and
\begin{eqnarray}\label{seteq10}
\frac{1}{n^2}\sum_{\substack{i\neq j\neq k,t\\ i_1\neq j_1\neq k_1,t_1\\ |S|=2}}\frac{1}{\nu_i^2\nu_{i_1}^2}\mathbb{E}\Bigg[B_{ijkt}B_{i_1j_1k_1t_1}\Bigg]
=O\left(\frac{n^5p_n^6}{n^2(np_n)^8}\right)
=O\left(\frac{1}{n^2(np_n)^3}\right).
\end{eqnarray}
If
$|S|=3$, then $|\{i,j,k,i_1,j_1,k_1\}|=3$  and
\begin{eqnarray}\label{seteq11}
\frac{1}{n^2}\sum_{\substack{i\neq j\neq k,t\\ i_1\neq j_1\neq k_1,t_1\\ |S|=3}}\frac{1}{\nu_i^2\nu_{i_1}^2}\mathbb{E}\Bigg[B_{ijkt}B_{i_1j_1k_1t_1}\Bigg]
=O\left(\frac{n^5p_n^5}{n^2(np_n)^8}\right)=O\left(\frac{1}{n^2(np_n)^3}\right).
\end{eqnarray}

Combining (\ref{seteq8})-(\ref{seteq11}) yields
\begin{eqnarray}\label{seteq12}
\mathbb{E}\left[\left(\frac{1}{n}\sum_{i\neq j\neq k,t}\frac{A_{ij}\bar{A}_{ij}A_{jk}A_{ki}\mu_{jt}-\mathbb{E}\left[A_{ij}\bar{A}_{ij}A_{jk}A_{ki}\mu_{jt}\right]}{\nu_i^2}\right)^2\right]=O\left(\frac{p_n}{n^2(np_n)^2}\right).
\end{eqnarray}

By (\ref{seteq6}), (\ref{seteq08}) and (\ref{seteq12}), we get
\begin{eqnarray}\nonumber
&&\frac{1}{n}\sum_{i\neq j\neq k, s\neq t}\frac{A_{ij}A_{jk}A_{ki}\bar{A}_{is}\mu_{st}-\mathbb{E}[A_{ij}A_{jk}A_{ki}\bar{A}_{is}\mu_{st}]}{\nu_i^2}\\ \label{seteq13}
&=&\frac{1}{n}\sum_{i\neq j\neq k\neq s,t}\frac{A_{ij}A_{jk}A_{ki}\bar{A}_{is}\mu_{st}}{\nu_i^2}+o_P\left(\frac{1}{n\sqrt{np_n}}+\frac{\sqrt{p_n}}{n}\right).
\end{eqnarray}

\vskip 1cm

Consider the last term of (\ref{seceq1}). Suppose $\{s,t\}\cap\{i,j,k\}=\emptyset$ in the summation of the last term of (\ref{seceq1}).
Similar to (\ref{seteq08}), we have
\begin{eqnarray}\label{seteq14}
\mathbb{E}\left[\left(\frac{1}{n}\sum_{i\neq j\neq k\neq s\neq t}\frac{A_{ij}A_{jk}A_{ki}\mu_{is}\bar{A}_{st}}{\nu_i^2}\right)^2\right]=O\left(\frac{p_n}{n^2}\right).
\end{eqnarray}
Suppose $\{s,t\}\subset\{i,j,k\}$ in the summation of the last term of (\ref{seceq1}). Since $s\neq i$ and $t\neq i$, then $\{s,t\}=\{j,k\}$. Hence
\begin{eqnarray}\label{seteq15}
\mathbb{E}\left[\left|\frac{1}{n}\sum_{i\neq j\neq k}\frac{A_{ij}A_{jk}A_{ki}\mu_{ij}\bar{A}_{jk}}{\nu_i^2}\right|\right]=O\left(\frac{1}{n^2}\right).
\end{eqnarray}

Suppose
$|\{s,t\}\cap\{i,j,k\}|=1$. Without loss of generality, let $s=j$ and $t\notin\{i,j,k\}$. By a similar argument as in (\ref{seteq3}), one has
\begin{eqnarray}\label{seteq16}
\mathbb{E}\left[\left(\frac{1}{n}\sum_{i\neq j\neq k\neq t}\frac{A_{ij}A_{jk}A_{ki}\mu_{ij}\bar{A}_{jt}}{\nu_i^2}\right)^2\right]=O\left(\frac{1}{n^2(np_n)^2}\right).
\end{eqnarray}

By (\ref{seteq14})-(\ref{seteq16}),
the last term of (\ref{seceq1}) is equal to
\begin{eqnarray}\label{seteq17}
\frac{1}{n}\sum_{i\neq j\neq k, s\neq t}\frac{A_{ij}A_{jk}A_{ki}\mu_{is}\bar{A}_{st}}{\nu_i^2}=\frac{1}{n}\sum_{i\neq j\neq k\neq s\neq t}\frac{A_{ij}A_{jk}A_{ki}\mu_{is}\bar{A}_{st}}{\nu_i^2}+o_P\left(\frac{\sqrt{p_n}}{n}\right).
\end{eqnarray}

Next we isolate the leading terms of  (\ref{seteq17}). It is easy to get that
\begin{eqnarray}\nonumber
\frac{1}{n}\sum_{i\neq j\neq k\neq s\neq t}\frac{A_{ij}A_{jk}A_{ki}\mu_{is}\bar{A}_{st}}{\nu_i^2}&=&\frac{1}{n}\sum_{i\neq j\neq k\neq s\neq t}\frac{(A_{ij}A_{jk}A_{ki}-\mu_{ij}\mu_{jk}\mu_{ki})\mu_{is}\bar{A}_{st}}{\nu_i^2}\\ \label{seteq18}
&&+\frac{1}{n}\sum_{i\neq j\neq k\neq s\neq t}\frac{\mu_{ij}\mu_{jk}\mu_{ki}\mu_{is}\bar{A}_{st}}{\nu_i^2}.
\end{eqnarray}

The second moment of the first term of (\ref{seteq18}) is equal to
\begin{eqnarray}\nonumber
&&\mathbb{E}\left[\left(\frac{1}{n}\sum_{i\neq j\neq k\neq s\neq t}\frac{(A_{ij}A_{jk}A_{ki}-\mu_{ij}\mu_{jk}\mu_{ki})\mu_{is}\bar{A}_{st}}{\nu_i^2}\right)^2\right]\\ \label{seteq19}
&=&\frac{1}{n^2}\sum_{\substack{i\neq j\neq k\neq s\neq t\\ i_1\neq j_1\neq k_1\neq s_1\neq t_1}}\frac{\mathbb{E}\left[B_{ijk}\mu_{is}\bar{A}_{st}B_{i_1j_1k_1}\mu_{i_1s_1}\bar{A}_{s_1t_1}\right]}{\nu_i^2\nu_{i_1}^2},
\end{eqnarray}
where $B_{ijk}=A_{ij}A_{jk}A_{ki}-\mu_{ij}\mu_{jk}\mu_{ki}$. Let 
\[S=\{(i,j),(j,k),(k,i),(s,t)\}\cap \{(i_1,j_1),(j_1,k_1),(k_1,i_1),(s_1,t_1)\}.\]

If $(s,t)\notin \{(i_1,j_1),(j_1,k_1),(k_1,i_1),(s_1,t_1)\}$, 
then 
\begin{equation*} 
\mathbb{E}\left[B_{ijk}\mu_{is}\bar{A}_{st}B_{i_1j_1k_1}\mu_{i_1s_1}\bar{A}_{s_1t_1}\right]=\mathbb{E}[\bar{A}_{st}]\mathbb{E}\left[B_{ijk}\mu_{is}B_{i_1j_1k_1}\mu_{i_1s_1}\bar{A}_{s_1t_1}\right]=0.
\end{equation*}
If
\[\{(i,j),(j,k),(k,i)\}\cap \{(i_1,j_1),(j_1,k_1),(k_1,i_1),(s_1,t_1)\}=\emptyset,\]
then 
\begin{equation}\label{seteq20}
\mathbb{E}\left[B_{ijk}\mu_{is}\bar{A}_{st}B_{i_1j_1k_1}\mu_{i_1s_1}\bar{A}_{s_1t_1}\right]=\mathbb{E}\left[B_{ijk}\right]\mathbb{E}\left[B_{ijk}\mu_{is}\bar{A}_{st}B_{i_1j_1k_1}\mu_{i_1s_1}\bar{A}_{s_1t_1}\right]=0.
\end{equation}

Hence
$(s,t)\in\{(i_1,j_1),(j_1,k_1),(k_1,i_1),(s_1,t_1)\}$ and
\[\{(i,j),(j,k),(k,i)\}\cap \{(i_1,j_1),(j_1,k_1),(k_1,i_1),(s_1,t_1)\}\neq\emptyset.\]
Then $|S|\geq2$. 
If $|S|=2$, then $|\{i,j,k,s,t,i_1,j_1,k_1,s_1,t_1\}|\leq 6$ and
\[\mathbb{E}\left[B_{ijk}\mu_{is}\bar{A}_{st}B_{i_1j_1k_1}\mu_{i_1s_1}\bar{A}_{s_1t_1}\right]=O(p_n^6).\]
Then
\begin{eqnarray}\label{seteq21}
\frac{1}{n^2}\sum_{\substack{i\neq j\neq k\neq s\neq t\\ i_1\neq j_1\neq k_1\neq s_1\neq t_1\\ |S|=2}}\frac{\mathbb{E}\left[B_{ijk}\mu_{is}\bar{A}_{st}B_{i_1j_1k_1}\mu_{i_1s_1}\bar{A}_{s_1t_1}\right]}{\nu_i^2\nu_{i_1}^2}=O\left(\frac{n^6p_n^6}{n^2(np_n)^8}\right)=O\left(\frac{1}{n^2(np_n)^2}\right).
\end{eqnarray}

If $|S|\geq3$, then $|\{i,j,k,s,t,i_1,j_1,k_1,s_1,t_1\}|\leq 5$ and
\[\mathbb{E}\left[B_{ijk}\mu_{is}\bar{A}_{st}B_{i_1j_1k_1}\mu_{i_1s_1}\bar{A}_{s_1t_1}\right]=O(p_n^5).\]
Then
\begin{eqnarray}\label{seteq22}
\frac{1}{n^2}\sum_{\substack{i\neq j\neq k\neq s\neq t\\ i_1\neq j_1\neq k_1\neq s_1\neq t_1\\ |S|\geq3}}\frac{\mathbb{E}\left[B_{ijk}\mu_{is}\bar{A}_{st}B_{i_1j_1k_1}\mu_{i_1s_1}\bar{A}_{s_1t_1}\right]}{\nu_i^2\nu_{i_1}^2}=O\left(\frac{n^5p_n^5}{n^2(np_n)^8}\right)=O\left(\frac{1}{n^2(np_n)^3}\right).
\end{eqnarray}

By (\ref{seteq19}), (\ref{seteq21}) and (\ref{seteq22}), we have
\begin{eqnarray}\label{seteq23}
\mathbb{E}\left[\left(\frac{1}{n}\sum_{i\neq j\neq k\neq s\neq t}\frac{(A_{ij}A_{jk}A_{ki}-\mu_{ij}\mu_{jk}\mu_{ki})\mu_{is}\bar{A}_{st}}{\nu_i^2}\right)^2\right]=O\left(\frac{1}{n^2(np_n)^2}\right).
\end{eqnarray}

Combining (\ref{seteq17}), (\ref{seteq18}) and (\ref{seteq23}) yields
\begin{eqnarray}\nonumber
&&\frac{1}{n}\sum_{i\neq j\neq k, s\neq t}\frac{A_{ij}A_{jk}A_{ki}\mu_{is}\bar{A}_{st}-\mathbb{E}[A_{ij}A_{jk}A_{ki}\mu_{is}\bar{A}_{st}]}{\nu_i^2}\\ \label{seteq24}
&=&\frac{1}{n}\sum_{i\neq j\neq k\neq s\neq t}\frac{\mu_{ij}\mu_{jk}\mu_{ki}\mu_{is}\bar{A}_{st}}{\nu_i^2}+o_P\left(\frac{1}{n\sqrt{np_n}}+\frac{\sqrt{p_n}}{n}\right).
\end{eqnarray}

Similarly, we have
\begin{eqnarray}\label{appeq}
\frac{1}{n}\sum_{i\neq j\neq k\neq s,t}\frac{A_{ij}A_{jk}A_{ki}\bar{A}_{is}\mu_{st}}{\nu_i^2}=\frac{1}{n}\sum_{i\neq j\neq k\neq s,t}\frac{\mu_{ij}\mu_{jk}\mu_{ki}\bar{A}_{is}\mu_{st}}{\nu_i^2}+o_P\left(\frac{1}{n\sqrt{np_n}}+\frac{\sqrt{p_n}}{n}\right).
\end{eqnarray}

\vskip 1cm

\noindent
{\bf Consider the third term of (\ref{eq1})}. 
Let $r_1,r_2,r_3$ be non-negative integers such that $r_1+r_2+r_3=t$. Let $C_{r_1r_2r_3}$ be constants dependent on $r_1,r_2,r_3$. By (\ref{lem2eq1}), we have
\begin{eqnarray}\label{thireq1}
    (V_i-\nu_i)^t=\sum_{r_1+r_2+r_3=t}C_{r_1r_2r_3}\left(\sum_{j\neq k}\bar{A}_{ij}\bar{A}_{jk}\right)^{r_1}\left(\sum_{j\neq k}\bar{A}_{ij}\mu_{jk}\right)^{r_2}\left(\sum_{j\neq k}\mu_{ij}\bar{A}_{jk}\right)^{r_3}.
\end{eqnarray}

It is easy to obtain that
\begin{eqnarray}\nonumber
A_{ij}A_{jk}A_{ki}
&=&\bar{A}_{ij}\bar{A}_{ik}\bar{A}_{jk}+\bar{A}_{ij}\bar{A}_{ik}\mu_{jk}+\bar{A}_{ij}\bar{A}_{jk}\mu_{ik}+\bar{A}_{ik}\bar{A}_{jk}\mu_{ij}\\ \label{neq15}
&&+\bar{A}_{jk}\mu_{ij}\mu_{ik}+\bar{A}_{ij}\mu_{jk}\mu_{ik}+\bar{A}_{ik}\mu_{jk}\mu_{ij}+\mu_{ij}\mu_{jk}\mu_{ki}.
\end{eqnarray}

Given non-negative integer $r\leq 2t$,
denote $D=\{i,j,k,j_1,\dots,j_t,k_1,\dots,k_t\}$ and
\[I_r=\Big\{(i,j,k,j_1,\dots,j_t,k_1,\dots,k_t)\in[n]^{2t+3}\Big| |D|=2t+3-r,i\neq j\neq k,j_s\neq k_s\neq i,s\in[t]\Big\}.\]
Given a set of indices in $D$, denote $e_1=(i,j)$, $e_2=(j,k)$, $e_3=(k,i)$, $e_{4}=(i,j_1)$, $e_5=(j_1,k_1)$, \dots, $e_{2t+2}=(i,j_t)$, $e_{2t+3}=(j_t,k_t)$ and
\[S=\big\{e_1,e_2,\dots,e_{2t+2},e_{2t+3}\},\]
where we assume $(i,j)=(j,i)$.
Then
\begin{eqnarray}\nonumber
&&\frac{1}{n}\sum_{i\neq j\neq k}^n\frac{1}{\nu_i^{t+1}}\bar{A}_{ij}\bar{A}_{ik}\bar{A}_{jk}\left(\sum_{j\neq k}\bar{A}_{ij}\bar{A}_{jk}\right)^{r_1}\left(\sum_{j\neq k}\bar{A}_{ij}\mu_{jk}\right)^{r_2}\left(\sum_{j\neq k}\mu_{ij}\bar{A}_{jk}\right)^{r_3}\\ \nonumber
&=& \frac{1}{n}\sum_{\substack{i\neq j\neq k\\j_1,\dots,j_t,\\ k_1,\dots,k_t}}\frac{\bar{A}_{ij}\bar{A}_{ik}\bar{A}_{jk}}{\nu_i^{t+1}}\left(\prod_{s=1}^{r_1}\bar{A}_{ij_s}\bar{A}_{j_sk_s}\right)\left(\prod_{s=r_1+1}^{r_1+r_2}\bar{A}_{ij_s}\mu_{j_sk_s}\right)\left(\prod_{s=r_1+r_2+1}^{t}\mu_{ij_s}\bar{A}_{j_sk_s}\right)\\ \nonumber
&=& \frac{1}{n}\sum_{I_0}\frac{\bar{A}_{ij}\bar{A}_{ik}\bar{A}_{jk}}{\nu_i^{t+1}}\left(\prod_{s=1}^{r_1}\bar{A}_{ij_s}\bar{A}_{j_sk_s}\right)\left(\prod_{s=r_1+1}^{r_1+r_2}\bar{A}_{ij_s}\mu_{j_sk_s}\right)\left(\prod_{s=r_1+r_2+1}^{t}\mu_{ij_s}\bar{A}_{j_sk_s}\right)\\ \nonumber
&&+ \frac{1}{n}\sum_{I_1}\frac{\bar{A}_{ij}\bar{A}_{ik}\bar{A}_{jk}}{\nu_i^{t+1}}\left(\prod_{s=1}^{r_1}\bar{A}_{ij_s}\bar{A}_{j_sk_s}\right)\left(\prod_{s=r_1+1}^{r_1+r_2}\bar{A}_{ij_s}\mu_{j_sk_s}\right)\left(\prod_{s=r_1+r_2+1}^{t}\mu_{ij_s}\bar{A}_{j_sk_s}\right)\\ \label{hardeq1}
&&+ \frac{1}{n}\sum_{I_r,r\geq2}\frac{\bar{A}_{ij}\bar{A}_{ik}\bar{A}_{jk}}{\nu_i^{t+1}}\left(\prod_{s=1}^{r_1}\bar{A}_{ij_s}\bar{A}_{j_sk_s}\right)\left(\prod_{s=r_1+1}^{r_1+r_2}\bar{A}_{ij_s}\mu_{j_sk_s}\right)\left(\prod_{s=r_1+r_2+1}^{t}\mu_{ij_s}\bar{A}_{j_sk_s}\right).
\end{eqnarray}
Next we bound each term in (\ref{hardeq1}).

Consider the last term in (\ref{hardeq1}). Suppose $r\geq2$. Let $e_{l},e_m\in S$ with $l\neq m$. If $e_l=e_m$, then there exist at least two elements $a,b\in D$ such that $a=b$. That is, when the number of distinct elements in $S$ reduces by 1, the number of distinct elements in $D$ reduces by at least 1. Hence $|S|<u$ implies $|D|<u$ for any integer $ 4\leq u\leq 2t+3$. Then $|D|\geq u$ implies $|S|\geq u$. Hence $|D|=2t+3-r$ implies $\big|S\big|\geq 2t+3-r$. In this case, one has
\[\mathbb{E}\left[\left|\bar{A}_{ij}\bar{A}_{ik}\bar{A}_{jk}\left(\prod_{s=1}^{r_1}\bar{A}_{ij_s}\bar{A}_{j_sk_s}\right)\left(\prod_{s=r_1+1}^{r_1+r_2}\bar{A}_{ij_s}\mu_{j_sk_s}\right)\left(\prod_{s=r_1+r_2+1}^{t}\mu_{ij_s}\bar{A}_{j_sk_s}\right)\right|\right]=O(p_n^{2t+3-r}).\]
Note that there are at most $n^{2t+3-r}$ elements in $I_r$. Then for $r\geq2$, we get
\begin{eqnarray}\nonumber
&&\frac{1}{n}\sum_{I_r,r\geq2}\mathbb{E}\left[\left|\frac{\bar{A}_{ij}\bar{A}_{ik}\bar{A}_{jk}}{\nu_i^{t+1}}\left(\prod_{s=1}^{r_1}\bar{A}_{ij_s}\bar{A}_{j_sk_s}\right)\left(\prod_{s=r_1+1}^{r_1+r_2}\bar{A}_{ij_s}\mu_{j_sk_s}\right)\left(\prod_{s=r_1+r_2+1}^{t}\mu_{ij_s}\bar{A}_{j_sk_s}\right)\right|\right]\\ \nonumber
&=&O\left(\frac{(np_n)^{2t+3-r}}{n(np_n)^{2t+2}}\right)\\ \label{hardeq2}
&=&O\left(\frac{1}{n(np_n)}\right).
\end{eqnarray}

\vskip 1cm

Consider the first term in (\ref{hardeq1}). 
Suppose $r=0$.
Let $D^{\prime}=\{i^{\prime},j^{\prime},k^{\prime},j_1^{\prime},\dots,j_t^{\prime},k_1^{\prime},\dots,k_t^{\prime}\}$ and
\[I_r^{\prime}=\Big\{(i^{\prime},j^{\prime},k^{\prime},j_1^{\prime},\dots,j_t^{\prime},k_1^{\prime},\dots,k_t^{\prime})\in[n]^{2t+3}\Big| |D^{\prime}|=2t+3-r,i^{\prime}\neq j^{\prime}\neq k^{\prime},j_s^{\prime}\neq k_s^{\prime}\neq i^{\prime},s\in [t]\Big\}.\]
Denote $e_1^{\prime}=(i^{\prime},j^{\prime})$, $e_2^{\prime}=(j^{\prime},k^{\prime})$, $e_3^{\prime}=(k^{\prime},i^{\prime})$, $e_{4}^{\prime}=(i^{\prime},j_1^{\prime})$, $e_5^{\prime}=(j_1^{\prime},k_1^{\prime})$, \dots, $e_{2t+2}^{\prime}=(i^{\prime},j_t^{\prime})$, $e_{2t+3}^{\prime}=(j_t^{\prime},k_t^{\prime})$, $S^{\prime}=\big\{e_1^{\prime},e_2^{\prime},\dots,e_{2t+2}^{\prime},e_{2t+3}^{\prime}\}$,
\[U=\big\{(j_{r_1+1},k_{r_1+1}),\dots,(j_{r_1+r_2},k_{r_1+r_2}), (i,j_{r_1+r_2+1}),\dots,(i,j_t)\},\]
\[U^{\prime}=\big\{(j_{r_1+1}^{\prime},k_{r_1+1}^{\prime}),\dots,(j_{r_1+r_2}^{\prime},k_{r_1+r_2}^{\prime}), (i^{\prime},j_{r_1+r_2+1}^{\prime}),\dots,(i^{\prime},j_t^{\prime})\}.\]
The second moment of the first term in (\ref{hardeq1}) is equal to
\begin{eqnarray}\nonumber
&&\mathbb{E}\left[\left(\frac{1}{n}\sum_{I_0}\frac{\bar{A}_{ij}\bar{A}_{ik}\bar{A}_{jk}}{\nu_i^{t+1}}\left(\prod_{s=1}^{r_1}\bar{A}_{ij_s}\bar{A}_{j_sk_s}\right)\left(\prod_{s=r_1+1}^{r_1+r_2}\bar{A}_{ij_s}\mu_{j_sk_s}\right)\left(\prod_{s=r_1+r_2+1}^{t}\mu_{ij_s}\bar{A}_{j_sk_s}\right)\right)^2\right]\\ \nonumber
&=& \frac{1}{n^2}\sum_{I_0,I_0^{\prime}}\mathbb{E}\Bigg[\frac{\bar{A}_{ij}\bar{A}_{ik}\bar{A}_{jk}}{\nu_i^{t+1}}\left(\prod_{s=1}^{r_1}\bar{A}_{ij_s}\bar{A}_{j_sk_s}\right)\left(\prod_{s=r_1+1}^{r_1+r_2}\bar{A}_{ij_s}\mu_{j_sk_s}\right)\left(\prod_{s=r_1+r_2+1}^{t}\mu_{ij_s}\bar{A}_{j_sk_s}\right)\\ \label{r1eq1}
&&\times \frac{\bar{A}_{i^{\prime}j^{\prime}}\bar{A}_{i^{\prime}k^{\prime}}\bar{A}_{j^{\prime}k^{\prime}}}{\nu_{i^{\prime}}^{t+1}}\left(\prod_{s=1}^{r_1}\bar{A}_{i^{\prime}j_s^{\prime}}\bar{A}_{j_s^{\prime}k_s^{\prime}}\right)\left(\prod_{s=r_1+1}^{r_1+r_2}\bar{A}_{i^{\prime}j_s^{\prime}}\mu_{j_s^{\prime}k_s^{\prime}}\right)\left(\prod_{s=r_1+r_2+1}^{t}\mu_{i^{\prime}j_s^{\prime}}\bar{A}_{j_s^{\prime}k_s^{\prime}}\right)\Bigg].
\end{eqnarray}

For any $e_l\in S-U$, $\bar{A}_{e_l}$ is independent of $\bar{A}_{e_m}$ for $e_m\in S-U-\{e_l\}$. If $\bar{A}_{e_l}\neq \bar{A}_{e_q^{\prime}}$ for any $e_q^{\prime}\in S^{\prime}-U^{\prime}$, then $\bar{A}_{e_l}$ is  independent of $\bar{A}_{e_m}$ for $e_m\in S-U-\{e_l\}$ and $\bar{A}_{e_q^{\prime}}$ for $e_q^{\prime}\in S^{\prime}-U^{\prime}$. In this case, the expectation in (\ref{r1eq1}) is zero. Hence, for the expectation in (\ref{r1eq1}) to be non-zero, we have that $\bar{A}_{e_l}=\bar{A}_{e_q^{\prime}}$ for some $e_q^{\prime}\in S^{\prime}-U^{\prime}$. That is, $e_l=e_q^{\prime}$ for some $e_q^{\prime}\in S^{\prime}-U^{\prime}$. Then we get 
\[S-U=S^{\prime}-U^{\prime}.\]
When $r=0$, the elements in $S-U$ are distinct. Hence each element in $S-U$ is equal to exactly one element in $S^{\prime}-U^{\prime}$.
Then by (\ref{r1eq1}) we have
\begin{eqnarray}\nonumber
&&\mathbb{E}\left[\left(\frac{1}{n}\sum_{I_0}\frac{\bar{A}_{ij}\bar{A}_{ik}\bar{A}_{jk}}{\nu_i^{t+1}}\left(\prod_{s=1}^{r_1}\bar{A}_{ij_s}\bar{A}_{j_sk_s}\right)\left(\prod_{s=r_1+1}^{r_1+r_2}\bar{A}_{ij_s}\mu_{j_sk_s}\right)\left(\prod_{s=r_1+r_2+1}^{t}\mu_{ij_s}\bar{A}_{j_sk_s}\right)\right)^2\right]\\ \nonumber
&=&O\left(\frac{1}{n^2(np_n)^{4t+4}}\right)\sum_{\substack{I_0,U^{\prime}\\ S-U=S^{\prime}-U^{\prime}}}\mathbb{E}\Bigg[\bar{A}_{ij}^2\bar{A}_{ik}^2\bar{A}_{jk}^2\left(\prod_{s=1}^{r_1}\bar{A}_{ij_s}^2\bar{A}_{j_sk_s}^2\right)\left(\prod_{s=r_1+1}^{r_1+r_2}\bar{A}_{ij_s}^2\mu_{j_sk_s}\right)\\ \label{thrr1eq2}
&&\times \left(\prod_{s=r_1+r_2+1}^{t}\mu_{ij_s}\bar{A}_{j_sk_s}^2\right)\left(\prod_{s=r_1+1}^{r_1+r_2}\mu_{j_s^{\prime}k_s^{\prime}}\right)\left(\prod_{s=r_1+r_2+1}^{t}\mu_{i^{\prime}j_s^{\prime}}\right)\Bigg].
\end{eqnarray}

Note that $S-U=S^{\prime}-U^{\prime}$ implies $j_s^{\prime}\in D$ for $s=r_1+1,\dots,t$ and $i^{\prime}\in D$. By (\ref{thrr1eq2}) and the fact that $r_2\leq t$, we have
\begin{eqnarray}\nonumber
&&\mathbb{E}\left[\left(\frac{1}{n}\sum_{I_0}\frac{\bar{A}_{ij}\bar{A}_{ik}\bar{A}_{jk}}{\nu_i^{t+1}}\left(\prod_{s=1}^{r_1}\bar{A}_{ij_s}\bar{A}_{j_sk_s}\right)\left(\prod_{s=r_1+1}^{r_1+r_2}\bar{A}_{ij_s}\mu_{j_sk_s}\right)\left(\prod_{s=r_1+r_2+1}^{t}\mu_{ij_s}\bar{A}_{j_sk_s}\right)\right)^2\right]\\ \nonumber
&=&O\left(\frac{n^{2t+3+r_2}p_n^{2t+3+r_2+r_3}}{n^2(np_n)^{4t+4}}\right)\\ \label{r1eq3}
&=&O\left(\frac{1}{n^2(np_n)^{t+1}}\right).
\end{eqnarray}

\vskip 1cm

Consider the second term in (\ref{hardeq1}). 
Suppose $r=1$. In this case, $ |D|=2t+3-1$. There are several cases: (a) $j_s=j$ or $j_s=k$ for some $s\in [r_1+r_2]$; (b) $j_s=j$ or $j_s=k$ for some $r_1+r_2+1\leq s\leq t$; (c) $j_{s_1}=j_{s_2}$ for distinct $s_1,s_2\in [r_1+r_2]$; (d) $j_{s_1}=j_{s_2}$ for $s_1\in [t]$ and $r_1+r_2+1\leq s_2\leq t$; (e)  $j_{s_1}=k_{s_2}$ for distinct $s_1,s_2$; (f) $k_{s_1}=k_{s_2}$ for distinct $s_1,s_2$; (g) $k_s=j$ or $k_s=k$.

Case (b). In this case, $|S|=2t+2$. Without loss of generality, we can assume $j_t=k$. In this case, $\bar{A}_{e_l}$ ($e_l\in S-U+\{i,k\}$) are independent and there are at most $n^{2t+2}$ elements in $I_1$.
By a similar argument as in (\ref{r1eq1})-(\ref{r1eq3}) yields
\begin{eqnarray}\nonumber
&&\mathbb{E}\left[\left(\frac{1}{n}\sum_{I_1,j_t=k}\frac{\bar{A}_{ij}\bar{A}_{ik}\bar{A}_{jk}}{\nu_i^{t+1}}\left(\prod_{s=1}^{r_1}\bar{A}_{ij_s}\bar{A}_{j_sk_s}\right)\left(\prod_{s=r_1+1}^{r_1+r_2}\bar{A}_{ij_s}\mu_{j_sk_s}\right)\left(\prod_{s=r_1+r_2+1}^{t}\mu_{ij_s}\bar{A}_{j_sk_s}\right)\right)^2\right]\\ \nonumber
&=&O\left(\frac{1}{n^2(np_n)^{4t+4}}\right)\sum_{\substack{I_1,j_t=k,U^{\prime}\\ S-U+\{i,k\}=S^{\prime}-U^{\prime}+\{i^{\prime},k^{\prime}\}}}\mathbb{E}\Bigg[\bar{A}_{ij}^2\bar{A}_{ik}^2\bar{A}_{jk}^2\left(\prod_{s=1}^{r_1}\bar{A}_{ij_s}^2\bar{A}_{j_sk_s}^2\right)\left(\prod_{s=r_1+1}^{r_1+r_2}\bar{A}_{ij_s}^2\mu_{j_sk_s}\right)\\  \nonumber\label{r1eq2}
&&\times \left(\prod_{s=r_1+r_2+1}^{t}\mu_{ij_s}\bar{A}_{j_sk_s}^2\right)\left(\prod_{s=r_1+1}^{r_1+r_2}\mu_{j_s^{\prime}k_s^{\prime}}\right)\left(\prod_{s=r_1+r_2+1}^{t}\mu_{i^{\prime}j_s^{\prime}}\right)\Bigg]\\ \nonumber
&=&O\left(\frac{n^{2t+2+r_2}p_n^{2t+3+r_2+r_3}}{n^2(np_n)^{4t+4}}\right)\\ \label{r2eq1}
&=&O\left(\frac{1}{n^2(np_n)^{t+2}}\right).
\end{eqnarray}
The cases (d), (e), (f), (g) can be similarly bounded as in (\ref{r2eq1}).

Case (a). Without loss of generality, we assume $j_1=j$. Then
\begin{eqnarray}\nonumber
&&\mathbb{E}\left[\left(\frac{1}{n}\sum_{I_1, j_1=j}\frac{\bar{A}_{ij}^2\bar{A}_{ik}\bar{A}_{jk}\bar{A}_{j_1k_1}}{\nu_i^{t+1}}\left(\prod_{s=2}^{r_1}\bar{A}_{ij_s}\bar{A}_{j_sk_s}\right)\left(\prod_{s=r_1+1}^{r_1+r_2}\bar{A}_{ij_s}\mu_{j_sk_s}\right)\left(\prod_{s=r_1+r_2+1}^{t}\mu_{ij_s}\bar{A}_{j_sk_s}\right)\right)^2\right]\\ \nonumber
&=&\frac{1}{n^2}\sum_{\substack{I_1,j_1=j\\ I_1^{\prime},j_1^{\prime}=j^{\prime}}}\mathbb{E}\Bigg[\frac{\bar{A}_{ij}^2\bar{A}_{ik}\bar{A}_{jk}\bar{A}_{j_1k_1}}{\nu_i^{t+1}}\left(\prod_{s=2}^{r_1}\bar{A}_{ij_s}\bar{A}_{j_sk_s}\right)\left(\prod_{s=r_1+1}^{r_1+r_2}\bar{A}_{ij_s}\mu_{j_sk_s}\right)\left(\prod_{s=r_1+r_2+1}^{t}\mu_{ij_s}\bar{A}_{j_sk_s}\right)\\ \label{r2eq2}
&&\times \frac{\bar{A}_{i^{\prime}j^{\prime}}^2\bar{A}_{i^{\prime}k^{\prime}}\bar{A}_{j^{\prime}k^{\prime}}\bar{A}_{j_1^{\prime}k_1^{\prime}}}{\nu_{i^{\prime}}^{t+1}}\left(\prod_{s=2}^{r_1}\bar{A}_{i^{\prime}j_s^{\prime}}\bar{A}_{j_s^{\prime}k_s^{\prime}}\right)\left(\prod_{s=r_1+1}^{r_1+r_2}\bar{A}_{i^{\prime}j_s^{\prime}}\mu_{j_s^{\prime}k_s^{\prime}}\right)\left(\prod_{s=r_1+r_2+1}^{t}\mu_{i^{\prime}j_s^{\prime}}\bar{A}_{j_s^{\prime}k_s^{\prime}}\right)\Bigg] 
\end{eqnarray}

For any $(l,m)\in S-U-\{(i,j)\}$, $\bar{A}_{lm}$ is independent of $\bar{A}_{ab}$ for $(a,b)\in S-U-\{(l,m)\}$. If $\bar{A}_{lm}\neq \bar{A}_{h^{\prime}g^{\prime}}$ for any $(h^{\prime},g^{\prime})\in S^{\prime}-U^{\prime}$, then $\bar{A}_{lm}$ is independent of $\bar{A}_{ab}$ for $(a,b)\in S-U-\{(l,m)\}$ and $(a,b)\in S^{\prime}-U^{\prime}$. In this case, the expectation in (\ref{r2eq2}) is zero. Hence, for the expectation in (\ref{r2eq2}) to be non-zero, we have that $\bar{A}_{lm}=\bar{A}_{h^{\prime}g^{\prime}}$ for some $(h^{\prime},g^{\prime})\in S^{\prime}-U^{\prime}$. That is, $(l,m)=(h^{\prime},g^{\prime})$ for some $(h^{\prime},g^{\prime})\in S^{\prime}-U^{\prime}$. Then we get 
\begin{equation}\label{r2eq3}
S-U-\{(i,j)\}\subset S^{\prime}-U^{\prime},\ \ \ \ \ \  S^{\prime}-U^{\prime}-\{(i^{\prime},j^{\prime})\}\subset S-U.
\end{equation}
There are at most $n^{2t+2+r_2}$ elements in $I_1$ and $I_1^{\prime}$ satisfying (\ref{r2eq3}). Note that $|\bar{A}_{ij}|\leq 1$. Then $|\bar{A}_{ij}|^m\leq\bar{A}_{ij}^2$ for $m\geq2$.   By (\ref{r2eq2}), we have
\begin{eqnarray}\nonumber
&&\mathbb{E}\left[\left(\frac{1}{n}\sum_{I_1,j_1=j}\frac{\bar{A}_{ij}^2\bar{A}_{ik}\bar{A}_{jk}\bar{A}_{j_1k_1}}{\nu_i^{t+1}}\left(\prod_{s=2}^{r_1}\bar{A}_{ij_s}\bar{A}_{j_sk_s}\right)\left(\prod_{s=r_1+1}^{r_1+r_2}\bar{A}_{ij_s}\mu_{j_sk_s}\right)\left(\prod_{s=r_1+r_2+1}^{t}\mu_{ij_s}\bar{A}_{j_sk_s}\right)\right)^2\right]\\ \nonumber
&\leq &\frac{1}{n^2}\sum_{\substack{I_1,j_1=j\\ I_1^{\prime},j_1^{\prime}=j^{\prime}\\ (\ref{r2eq3})}}\mathbb{E}\Bigg[\frac{\bar{A}_{ij}^2\bar{A}_{ik}^2\bar{A}_{jk}^2\bar{A}_{j_1k_1}^2}{\nu_i^{t+1}\nu_{i^{\prime}}^{t+1}}\left(\prod_{s=2}^{r_1}\bar{A}_{ij_s}^2\bar{A}_{j_sk_s}^2\right)\left(\prod_{s=r_1+1}^{r_1+r_2}\bar{A}_{ij_s}^2\mu_{j_sk_s}\right)\left(\prod_{s=r_1+r_2+1}^{t}\mu_{ij_s}\bar{A}_{j_sk_s}^2\right)\\     \nonumber
&&\times\left(\prod_{s=r_1+1}^{r_1+r_2}\mu_{j_s^{\prime}k_s^{\prime}}\right)\left(\prod_{s=r_1+r_2+1}^{t}\mu_{i^{\prime}j_s^{\prime}}\right)\Bigg]\\  \nonumber
&=&O\left(\frac{n^{2t+2+r_2}p_n^{2t+2+r_2+r_3}}{n^2(np_n)^{4t+4}}\right)\\  \label{r2eq4}
&=&O\left(\frac{1}{n^2(np_n)^{t+2}}\right).
\end{eqnarray}

Case (c) can be similarly bounded as in (\ref{r2eq4}).

By (\ref{hardeq1}), (\ref{hardeq2}), (\ref{r1eq3}), (\ref{r2eq1}) and (\ref{r2eq4}), we get
\begin{eqnarray}\nonumber
&&\frac{1}{n}\sum_{i\neq j\neq k}^n\frac{\bar{A}_{ij}\bar{A}_{ik}\bar{A}_{jk}}{\nu_i^{t+1}}\left(\sum_{j\neq k\neq i}\bar{A}_{ij}\bar{A}_{jk}\right)^{r_1}\left(\sum_{j\neq k\neq i}\bar{A}_{ij}\mu_{jk}\right)^{r_2}\left(\sum_{j\neq k\neq i}\mu_{ij}\bar{A}_{jk}\right)^{r_3}\\ \nonumber
&&-\mathbb{E}\left[\frac{1}{n}\sum_{i\neq j\neq k}^n\frac{\bar{A}_{ij}\bar{A}_{ik}\bar{A}_{jk}}{\nu_i^{t+1}}\left(\sum_{j\neq k\neq i}\bar{A}_{ij}\bar{A}_{jk}\right)^{r_1}\left(\sum_{j\neq k\neq i}\bar{A}_{ij}\mu_{jk}\right)^{r_2}\left(\sum_{j\neq k\neq i}\mu_{ij}\bar{A}_{jk}\right)^{r_3}\right] \\ \label{hardeq5}
&=&o_P\left(\frac{1}{n\sqrt{np}}\right). 
\end{eqnarray}
Similarly, it is easy to show that
\begin{eqnarray}\nonumber
&&\frac{1}{n}\sum_{i\neq j\neq k}^n\frac{B_{ijk}}{\nu_i^{t+1}}\left(\sum_{j\neq k\neq i}\bar{A}_{ij}\bar{A}_{jk}\right)^{r_1}\left(\sum_{j\neq k\neq i}\bar{A}_{ij}\mu_{jk}\right)^{r_2}\left(\sum_{j\neq k\neq i}\mu_{ij}\bar{A}_{jk}\right)^{r_3}\\ \nonumber
&&-\mathbb{E}\left[\frac{1}{n}\sum_{i\neq j\neq k}^n\frac{B_{ijk}}{\nu_i^{t+1}}\left(\sum_{j\neq k\neq i}\bar{A}_{ij}\bar{A}_{jk}\right)^{r_1}\left(\sum_{j\neq k\neq i}\bar{A}_{ij}\mu_{jk}\right)^{r_2}\left(\sum_{j\neq k\neq i}\mu_{ij}\bar{A}_{jk}\right)^{r_3}\right]\\ \label{hardeq6}
&=&o_P\left(\frac{1}{n\sqrt{np}}\right),
\end{eqnarray}
where $B_{ijk}\in\{\bar{A}_{ij}\bar{A}_{ik}\mu_{jk},\bar{A}_{ij}\mu_{ik}\bar{A}_{jk},\mu_{ij}\bar{A}_{ik}\bar{A}_{jk},\bar{A}_{ij}\mu_{ik}\mu_{jk},\mu_{ij}\bar{A}_{ik}\mu_{jk},\mu_{ij}\mu_{ik}\bar{A}_{jk},\mu_{ij}\mu_{ik}\mu_{jk}\}$.

Based on (\ref{neq15}), (\ref{thireq1}), (\ref{hardeq5}) and (\ref{hardeq6}), the third term of (\ref{eq1}) is bounded by
\begin{eqnarray}
\sum_{t=2}^{k_0}\frac{1}{n}\sum_{i=1}^n\frac{\Delta_{i}(V_i-\nu_i)^t-\mathbb{E}\left[\Delta_{i}(V_i-\nu_i)^t\right]}{\nu_i^{t+1}}=o_P\left(\frac{1}{n\sqrt{np}}\right).
\end{eqnarray}

By (\ref{y1ng}), (\ref{y1ns}), (\ref{y1eq}),  (\ref{seteq24}) and (\ref{appeq}), we have the following result.

If $\alpha>\frac{1}{2}$, then
\begin{eqnarray}\label{y1ng}
\overline{\mathcal{H}}_n-\mathbb{E}[\overline{\mathcal{H}}_n]
&=&\frac{2}{n}\sum_{i<j<k}\left(\frac{1}{\nu_i}+\frac{1}{\nu_j}+\frac{1}{\nu_k}\right)\bar{A}_{ij}\bar{A}_{ik}\bar{A}_{jk}+o_P\left(\frac{1}{n\sqrt{np_n}}\right).
\end{eqnarray}

If $\alpha<\frac{1}{2}$, then
\begin{eqnarray} \label{y1ns}
\overline{\mathcal{H}}_n-\mathbb{E}[\overline{\mathcal{H}}_n]
=\frac{1}{n}\sum_{i< j}(2b_{ij}+2c_{ij}+2c_{ji}-(a_i+a_j)-(e_{ij}+e_{ji}))\bar{A}_{ij}+o_P\left(\frac{\sqrt{p_n}}{n}\right).
\end{eqnarray}

If $\alpha=\frac{1}{2}$, then
\begin{eqnarray} \nonumber
&&\overline{\mathcal{H}}_n-\mathbb{E}[\overline{\mathcal{H}}_n]\\
&=&\frac{2}{n}\sum_{i<j<k}\left(\frac{1}{\nu_i}+\frac{1}{\nu_j}+\frac{1}{\nu_k}\right)\bar{A}_{ij}\bar{A}_{ik}\bar{A}_{jk}\\
&&+\frac{1}{n}\sum_{i< j}(2b_{ij}+c_{ij}+c_{ji}-(a_i+a_j)-(e_{ij}+e_{ji}))\bar{A}_{ij}+o_P\left(\frac{\sqrt{p_n}}{n}\right).
\end{eqnarray}

By a similar proof of Lemma 3.4 in \cite{YZ23}, we get the limiting distribution of $\overline{\mathcal{H}}_n-\mathbb{E}[\overline{\mathcal{H}}_n]$ as given in Theorem \ref{mainthm}.

\qed

\end{document}